\documentclass{article}

\usepackage{epsfig}
\usepackage{graphicx}
\usepackage{amsmath}
\usepackage{amsfonts}
\usepackage{amssymb}
\usepackage{tikz}
\usepackage{amsthm}
\usepackage{subfigure}

\newtheorem{thm}{Theorem}[section]

\newtheorem{rek}[thm]{Remark}

\DeclareMathOperator*{\argmin}{arg\,min}

\title{A SYMBOL-BASED ALGORITHM FOR DECODING BAR
  CODES}  

\author{Mark A.~Iwen\thanks{Mathematics Department, Duke University, 
Durham, NC 27708 (markiwen@math.duke.edu).  The research of this author was
supported in part by ONR N00014-07-1-0625 and NSF DMS
DMS-0847388.}  
\and Fadil Santosa\footnotemark\thanks{School of
Mathematics, University of Minnesota, Minneapolis, MN 55455
(santosa@math.umn.edu).  The research of this author was supported
in part by NSF DMS-0807856.}  
\and Rachel Ward\thanks{Department
of Mathematics, University of Texas at Austin, 2515 Speedway,
Austin, TX, 78712 (rward@math.utexas.edu).  The research of this
author was supported in part by the NSF Postdoctoral Research
Fellowship and the Donald D. Harrington Faculty Fellowship.}}

\begin{document}
\maketitle

\pagestyle{myheadings}â
\thispagestyle{plain}â
\markboth{M.~A.~IWEN, F.~SANTOSA, AND R.~WARD}{SYMBOL-BASED BAR CODE
  DECODING}

\begin{abstract}
We investigate the problem of decoding a bar code from a signal
measured with a hand-held laser-based scanner.  Rather than
formulating the inverse problem as one of binary image reconstruction,
we instead incorporate the symbology of the bar code into the
reconstruction algorithm directly, and search for a sparse
representation of the UPC bar code with respect to this known
dictionary.  Our approach significantly reduces the degrees of freedom
in the problem, allowing for accurate reconstruction that is robust to
noise and unknown parameters in the scanning device.  We propose a
greedy reconstruction algorithm and provide robust reconstruction
guarantees.  Numerical examples illustrate the insensitivity of our
symbology-based reconstruction to both imprecise model parameters and
noise on the scanned measurements.
\end{abstract}



\section{Introduction}

This work concerns an approach for decoding bar code signals.  While
it is true that bar code scanning is essentially a solved problem in
many domains, as evidenced by its prevalent use, there is still a need
for more reliable decoding algorithms in situations where the signals
are highly corrupted and the scanning takes place in less than ideal
situations.  It is under these conditions that traditional bar code
scanning algorithms often fail.

The problem of bar code decoding may be viewed as the deconvolution of
a binary one-dimensional image involving unknown parameters in the
blurring kernel that must be estimated from the signal
\cite{esedoglu}.  Esedoglu \cite{esedoglu} was the first to provide a
mathematical analysis of the bar code decoding problem in this
context, and he established the first uniqueness result of its kind
for this problem.  He further showed that the blind deconvolution
problem can be formulated as a well-posed variational problem.  An
approximation, based on the Modica-Mortola energy
\cite{modica-mortola}, is the basis for the computational approach.
The approach has recently been given further analytical treatment in
\cite{esedoglu-santosa}.

A recent work \cite{choksi-vangennip} addresses the case where the
blurring is not very severe.  Indeed the authors were able to treat
the signal as if it has not been blurred.  They showed rigorously
the variational framework
can recover the true bar code image even if this parameter is not
known.  A later paper \cite{choksi-vangennip-oberman} consider the case
where blurring is large and its parameter value known.  However, none of these
papers deal rigorously with noise although their numerical simulations
included noise.  For an analysis of the deblurring problem where the
blur is large and noise is present, the reader is referred to
\cite{esedoglu-santosa}. 

The approach presented in this work departs from the above image-based
approaches.  We treat the unknown as a finite-dimensional code and
develop a model that relates the code to the measured signal.  We show
that by exploiting the symbology -- the language of the bar code -- a
bar code can be identified with a sparse representation in the
symbology dictionary.  We develop a recovery algorithm that fits the
observed signal to a code from the symbology in a greedy fashion,
iterating in one pass from left to right.  We prove that the algorithm
can tolerate a significant level of blur and noise.  We also verify
insensitivity of the reconstruction to imprecise parameter estimation
of the blurring function.

We were unable to find any previous symbol-based methods for bar code
decoding in the open literature.  In a related approach
\cite{dumas-Rhabi-Rochefort}, a genetic algorithm is utilized to
represent populations of candidate barcodes together with likely
blurring and illumination parameters from the observed image data.
Successive generations of candidate solutions are then spawned from
those best matching the input data until a stopping criterion is met.
That work differs from the current article in that it uses a different
decoding method and does not utilize the relationship between the
structure of the barcode symbology and the blurring kernel. 

We note that there is a symbol-based approach for super-resolving
scanned images \cite{bern-goldberg}.  However, that work is
statistical in nature whereas the method we present is deterministic.
Both this work and the super-resolution work are similar in spirit to
lossless data compression algorithms known as `dictionary coding'
(see, e.g., \cite{salomon2004data}) which involve matching strings of
text to strings contained in an encoding dictionary.

The outline of the paper is as follows.  We start by developing a
model for the scanning process.  In Section 3, we study the
properties of the UPC (Universal Product Code) 
bar code and provide a mathematical representation
for the code. Section 4 develops the relation between the code and
the measured signal. An algorithm for decoding bar code signals
is presented in Section 5.  Section 6 is devoted to the analysis
of the algorithm proposed.  Results from numerical experiments
are presented in Section 7, and a final section concludes the work
with a discussion.

\section{A scanning model and associated inverse problem}

A bar code is scanned by shining a narrow laser beam 
across the black-and-white bars at constant speed.  The amount of light reflected
as the beam moves is recorded and can be viewed as a signal in
time. Since the bar code consists of black and white segments, the
reflected energy is large when the beam is on the white part, and
small when the beam is on the black part.  The reflected light energy 
at a given position is proportional to the integral
of the product of the beam intensity, which can be modeled as a Gaussian\footnote{This Gaussian model has also been utilized in many previous treatments of the bar code decoding problem.  See, e.g., \cite{kim-lee} and references therein.}, 
and the bar code image intensity
(white is high intensity, black is low).
The recorded data are samples of the resulting
continuous time signal.

\begin{figure}[ht]
\centerline{\includegraphics[width=5in]{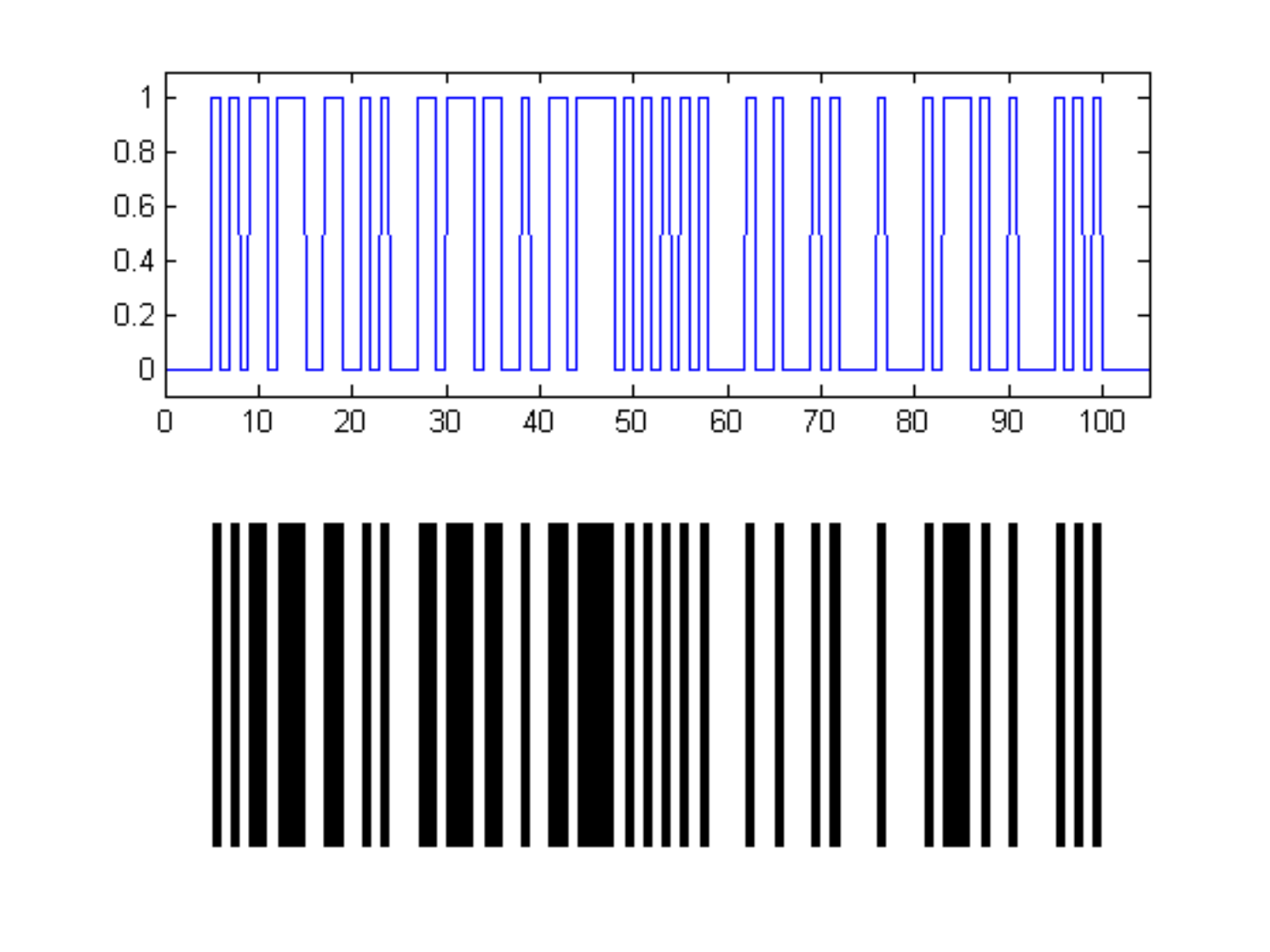}}
\caption{Samples of the binary bar code function $z(t)$ and the UPC bar code.
Note that in UPC bar codes, each bar - black or white - is a multiple of
the minimum bar width.}
\label{fig:ExampBar}
\end{figure}

Let us write the Gaussian beam intensity as a function of time:
\begin{equation}\label{thegaussian}
g(t) = \alpha
\frac{1}{\sqrt{2\pi}\sigma} e^{-(t^2/2\sigma^2)}.
\end{equation}
There are two parameters: (i) the variance $\sigma^2$ and (ii) the
constant multiplier $\alpha$.    We
will overlook the issue of relating time to the actual position of the
laser beam on the bar code, which is measured in distance. We can do
this because only relative widths of the bars are important in their encoding.  

Because the bar code -- denoted by $z(t)$ -- represents a black and white image,
we will normalize it to be a binary function.
Then the sampled data are
\begin{equation}\label{dataeqn}
d_i = \int g(t_i -\tau) z(\tau) d\tau + h_i , \hspace{3mm} i \in [m],
\end{equation}
where the $t_i \in [0,n]$ are equally spaced discretization points, and the $h_i$ represent the noise associated with scanning.  We have used the 
notation $[m] = \{1,2, ..., m\}$.  We need to
consider the relative size of the laser beam spot to the width of the
narrowest bar in the bar code. We set the minimum bar width to be 1
in the artificial time measure. 

It remains to explain the roles of the parameters in the Gaussian beam intensity. The
variance $\sigma^2$ models the distance from the scanner to the bar code, with larger variance
signifying longer distance.  The
\emph{width} of a Gaussian represents the length of the interval, centered around the Gaussian mean, over which the
Gaussian is greater than half its maximum amplitude; it is given by $2\sqrt{2\ln 2} \sigma$.
Informally, the Gaussian blur width should be of the same order of magnitude as
the size as the minimum bar width in the bar code for possible reconstruction. The multiplier $\alpha$ lumps the
conversion from light energy interacting with a binary bar code image to
the measurement.  Since the distance to the bar code is unknown and
the intensity-to-voltage conversion depends on ambient light and
properties of the laser/detector, these parameters are assumed to
be unknown.

To develop the model further, consider
the characteristic function
\[
\chi(t) = \left\{ \begin{array}{ll}
1 & \mbox{for} \;\; 0 \leq t \leq 1, \\
0 & \mbox{else}.
                  \end{array} \right.
\]
Then the bar code function can be written as
\begin{equation}\label{barcoderep}
z(t) = \sum_{j=1}^n c_j \chi(t - (j-1)),
\end{equation}
where the coefficients $c_j$ are either 0 or 1 (see, e.g., Figure~\ref{fig:ExampBar}).  The sequence
\[
c_1, c_2, \cdots, c_n,
\]
represents the information stored in the bar code, with a `0' corresponding to a white bar of unit width 
and a `1' corresponding to a black bar of unit width. 
For UPC bar codes, the total number of unit widths, $n$, is fixed to be $95$ 
for a $12$-digit code (further
explanations in the subsequent).

\begin{rek} \textnormal{One can think of the sequence $\{c_1,c_2,\cdots,c_n\}$ as
  an instruction for printing a bar code. Every $c_i$ is a command to
  lay out a white bar if $c_i=0$, or a black bar if otherwise.}
\end{rek}

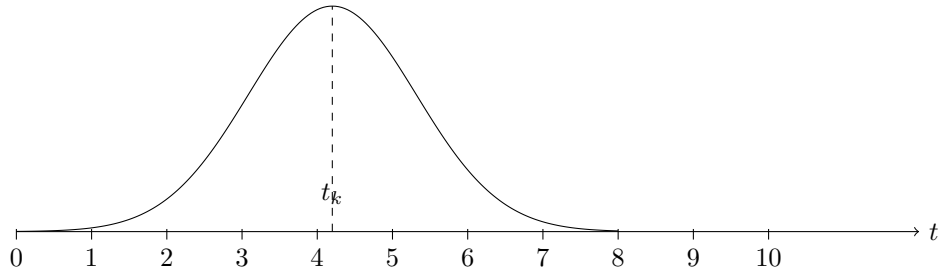
\begin{figure}[ht]
\label{fig:gkj}
\begin{center}
\begin{tikzpicture}
\draw[smooth,samples=100,domain=-4:4]
plot(\x,{3*exp(-0.4*(\x-0.2)^2) });
\draw[dashed] (0.2,0) -- (0.2,3);
\draw[->] (-4,0) -- (8,0);
\draw (0.2,0.5) node{$t_k$};
\foreach \x in {0,...,10}
     		\draw (\x-4,1pt) -- (\x-4,-3pt)
			node[anchor=north] {\x};
\draw (8.2,0) node{$t$};
\end{tikzpicture}
\end{center}
\caption{The matrix element ${\cal G}_{kj}$ is calculated by
    placing a scaled Gaussian over the bar code grid and integrating
    over each of the bar code intervals.}
\end{figure}

\noindent Substituting the bar code representation \eqref{barcoderep}
back in (\ref{dataeqn}), the sampled data can be represented as follows:
\begin{align*}
d_i &= \int g(t_i-t) \left[ \sum_{j=1}^n c_j \chi(t - (j-1)) \right]
dt + h_i \\
&= \sum_{j=1}^n \left[ \int_{(j-1)}^j g(t_i-t) dt \right] c_j  + h_i .
\end{align*}
In terms of the matrix ${\cal G} = {\cal G}(\sigma)$ with entries 
\begin{equation}\label{defineG}
{\cal G}_{kj} = \frac{1}{\sqrt{2\pi}\sigma}\int_{(j-1)}^j e^{-\frac{(t_k-t)^2}{2\sigma^2}}dt, \;\; \;\;k \in [m], \;\; j \in [n],
\end{equation}
the bar code determination problem reads
\begin{equation}\label{inverseprob}
d = \alpha {\cal G}(\sigma) c + h.
\end{equation}
The matrix entries ${\cal G}_{kj}$ are illustrated in Figure~2.\ref{fig:gkj}.
In the sequel, we will assume this discrete version of the bar code problem.
While it is tempting to solve (\ref{inverseprob}) directly for $c$,
$\sigma$ and $\alpha$, the best approach for doing so is not obvious.
The main difficulty stems from the fact that $c$ is a binary vector,
while the Gaussian parameters are continuous variables.  

\section{Incorporating the UPC bar code symbology}

We now tailor the bar code reading problem to UPC bar codes, although we remark that our general
framework should apply generally to any bar code of fixed length.
In the UPC-A symbology, a bar code represents a
12-digit number.  If we ignore the check-sum requirement, then any
12-digit number is permitted, and the number of unit widths, $n$, is fixed to $95$.  Going from left to right, the UPC bar code
has 5 parts -- the start sequence, the codes for the first 6 digits,
the middle sequence, the codes for the next 6 digits, and the end
sequence.  Thus the bar code has the following structure:
\begin{equation}
SL_1L_2L_3L_4L_5L_6MR_1R_2R_3R_4R_5R_6E,
\label{barcodestruct}
\end{equation}
where $S$, $M$, and $E$ are the start, middle, and end patterns
respectively, and
$L_i$ and $R_i$ are patterns corresponding to the digits.

In the sequel, we represent a white bar of unit width by 0 and a
black bar by 1 in the bar code representation $\{c_i\}$.\footnote{Note that identifying white bars with 0 and black bars with 1 runs 
counter to the natural light intensity of the reflected laser beam.  However, it is the black bars that carry information.}
The start, middle, and end patterns are fixed and given by
\[
S=E=[101], \;\; M=[01010].
\]
The patterns for $L_i$ and $R_i$ are taken from the following table:
\begin{equation}
\begin{tabular}{|c|r|r|}\hline
digit & L-pattern & R-pattern \\ \hline \hline
0     & $0001101$ & $1110010$ \\ \hline
1     & $0011001$ & $1100110$ \\ \hline
2     & $0010011$ & $1101100$ \\ \hline
3     & $0111101$ & $1000010$ \\ \hline
4     & $0100011$ & $1011100$ \\ \hline
5     & $0110001$ & $1001110$ \\ \hline
6     & $0101111$ & $1010000$ \\ \hline
7     & $0111011$ & $1000100$ \\ \hline
8     & $0110111$ & $1001000$ \\ \hline
9     & $0001011$ & $1110100$ \\ \hline
\end{tabular}
\label{LRcodes}
\end{equation}
Note that the right patterns are just the left patterns with the
0's and 1's flipped.  It follows that the bar
code can be represented as a binary vector $c \in \{0,1\}^{95}$.
However, not every binary vector constitutes a bar code -- only 
$10^{12}$ of the possible $2^{95}$ binary sequences of length $95$ -- fewer than $10^{-16} \hspace{.5mm} \%$ -- are bar codes. Specifically, the bar code structure (\ref{barcodestruct}) indicates 
that bar codes have specific sparse representations in the bar code dictionary constructed as
follows: write the
left-integer and right-integer codes as columns of a 7-by-10 matrix,
\[
L = \left[ \begin{array}{cccccccccc}
0 & 0 & 0 & 0 & 0 & 0 & 0 & 0 & 0 & 0 \\
0 & 0 & 0 & 1 & 1 & 1 & 1 & 1 & 1 & 0 \\
0 & 1 & 1 & 1 & 0 & 1 & 0 & 1 & 1 & 0 \\
1 & 1 & 0 & 1 & 0 & 0 & 1 & 1 & 0 & 1 \\
1 & 0 & 0 & 1 & 0 & 0 & 1 & 0 & 1 & 0 \\
0 & 0 & 1 & 0 & 1 & 0 & 1 & 1 & 1 & 1 \\
1 & 1 & 1 & 1 & 1 & 1 & 1 & 1 & 1 & 1 
            \end{array}  \right]  ,
\]

\[
R = \left[ \begin{array}{cccccccccc}
1 & 1 & 1 & 1 & 1 & 1 & 1 & 1 & 1 & 1 \\
1 & 1 & 1 & 0 & 0 & 0 & 0 & 0 & 0 & 1 \\
1 & 0 & 0 & 0 & 1 & 0 & 1 & 0 & 0 & 1 \\
0 & 0 & 1 & 0 & 1 & 1 & 0 & 0 & 1 & 0 \\
0 & 1 & 1 & 0 & 1 & 1 & 0 & 1 & 0 & 1 \\
1 & 1 & 0 & 1 & 0 & 1 & 0 & 0 & 0 & 0 \\
0 & 0 & 0 & 0 & 0 & 0 & 0 & 0 & 0 & 0 
            \end{array}  \right]  .
\]
The start and end patterns, $S$ and $E$, are 3-dimensional vectors,
while the middle pattern $M$ is a 5-dimensional vector 
\[
S = E= [0 1 0]^T, \;\; M = [ 0 1 0 1 0]^T.
\]
The bar code dictionary is the 95-by-123 block diagonal matrix
\[
{\cal D} = 
\left[ \begin{array}{ccccccccccccccc}
S  & 0 & \dots&    &   &   &   &   &   &   &   &   &   &\dots& 0  \\
0  & L &   &   &   &   &   &   &   &   &   &   &   &   &\vdots  \\
\vdots &   & L &   &   &   &   &   &   &   &   &   &   &   &  \\
   &   &   & L &   &   &   &   &   &   &   &   &   &   &  \\
   &   &   &   & L &   &   &   &   &   &   &   &   &   &  \\
   &   &   &   &   & L &   &   &   &   &   &   &   &   &  \\
   &   &   &   &   &   & L &   &   &   &   &   &   &   &  \\
   &   &   &   &   &   &   & M &   &   &   &   &   &   &  \\
   &   &   &   &   &   &   &   & R &   &   &   &   &   &  \\
   &   &   &   &   &   &   &   &   & R &   &   &   &   &  \\
   &   &   &   &   &   &   &   &   &   & R &   &   &   &  \\
   &   &   &   &   &   &   &   &   &   &   & R &   &   &  \\
   &   &   &   &   &   &   &   &   &   &   &   & R &   &\vdots\\
\vdots &   &   &   &   &   &   &   &   &   &   &   &   & R & 0 \\
0   &\dots   &   &   &   &   &   &   &   &   &   &   &\dots & 0  & E 
\end{array} \right] .
\]
The bar code (\ref{barcodestruct}), expanded in the bar code dictionary, has the form
\begin{equation}\label{dictrep}
c = {\cal D} x, \quad \quad x \in \{0,1\}^{123},
\end{equation}
where 
\begin{enumerate}{\itemsep=0cm\parsep=0cm}
\item The 1st, 62nd and the 123rd entries of $x$, corresponding to the
  $S$, $M$, and $E$ patterns, are $1$.
\item Among the 2nd through 11th entries of $x$, exactly one
  entry -- the entry corresponding to the first digit in $c = {\cal D}x$ -- is
  nonzero. The same is true for 12th through 22nd entries, etc, until
  the 61st entry.  This pattern starts again from the 63rd entry
  through the 122th entry. In all, $x$ has exactly $15$ nonzero entries.
\end{enumerate}
That is, $x$ must take the form
\begin{equation}\label{xform}
x^T = [1, v_1^T, \cdots, v_6^T, 1, v_7^T, \cdots, v_{12}^T, 1],
\end{equation}
where $v_j$, for $j=1,\cdots,12$, are vectors in $\{0,1\}^{10}$
having only one nonzero element.
\noindent In this new representation,
the bar code reconstruction problem \eqref{inverseprob} reads
\begin{equation}
d = \alpha {\cal G}(\sigma){\cal D} x + h,
\label{noisydata}
\end{equation}
where $d \in \mathbb{R}^m$ is the measurement vector, the matrices 
${\cal G}\textrm{(}\sigma \textrm{)} \in \mathbb{R}^{m \times 95} $ 
and ${\cal D} \in \{ 0,1 \}^{95 \times 123}$ are as defined in (\ref{defineG}) 
and (\ref{dictrep}) respectively, and 
$h \in \mathbb{R}^m$ is additive noise. 
Note that ${\cal D}$ has fewer rows than columns, while ${\cal G}$ will generally 
have more rows than columns; 
we will refer to the ratio of rows to columns as the oversampling ratio and denote 
it by $r = m/n$.  Given the data
$d \in \mathbb{R}^{m}$, our objective is to return a valid bar code 
 $x \in \{ 0, 1 \}^{123}$ as reliably and 
quickly as possible.  

\section{Properties of the forward map}

Incorporating the bar code dictionary into the inverse problem \eqref{noisydata}, we see that the map
between the bar code and observed data is represented by the matrix
${\cal P} = \alpha {\cal G}( \sigma) {\cal D} \in \mathbb{R}^{m \times 123}$.  We will
refer to ${\cal P}$, which is a function of the model parameters $\alpha$ and $\sigma$, 
as the \emph{forward map}. 

\subsection{Near block-diagonality} 
\begin{figure}[h]
\centerline{
\includegraphics[width=2.5in]{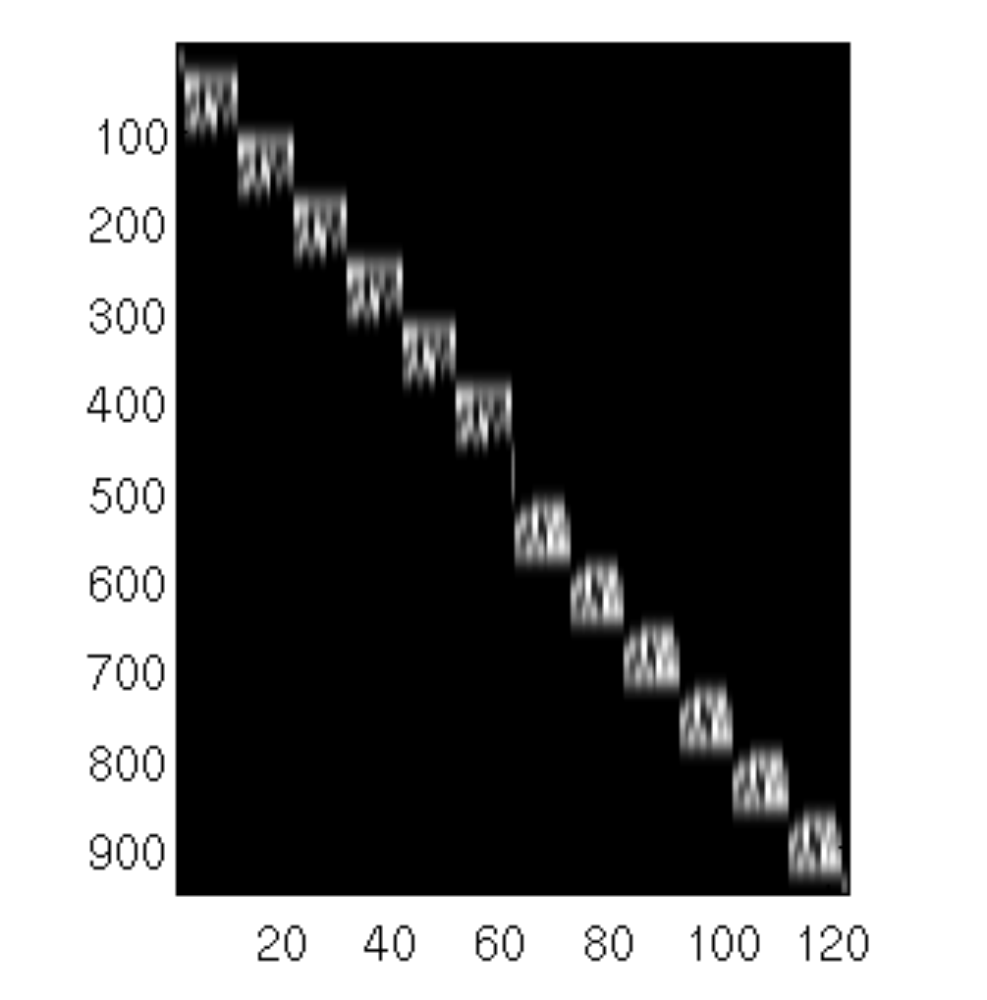} }
\caption{A representative bar code forward map ${\cal P} = \alpha{\cal
    G}(\sigma){\cal D}$ corresponding to oversampling parameter $r =
  10$, amplitude $\alpha = 1$, and Gaussian standard deviation $\sigma
  = 1.5$.  The lone column vectors at the start, middle, and end
  account for the known start, middle, and end patterns in the bar
  code.}
\label{fig:3}
\end{figure}

For reasonable levels of blur in the Gaussian kernel,
the forward map ${\cal P}$ inherits an almost block-diagonal structure
from the bar code matrix ${\cal D}$ as illustrated in Figure \ref{fig:3}.
 In the limit as the amount of blur $\sigma \rightarrow 0$, the forward map
${\cal P}$ becomes exactly the block-diagonal bar code matrix. 
More precisely, we partition the forward map ${\cal P}$ according to 
the block-diagonal structure of the bar code dictionary ${\cal D}$,

\begin{equation} \label{defineP}
{\cal P} = \left[ P^{(1)} ~P^{(2)} ~\dots ~P^{(15)} \right] .
\end{equation}
The 1st, 8th, and 15th sub-matrices are special as they correspond to
the known start, middle, and end patterns of the bar code.  In accordance with the 
structure of $x$ where $c = {\cal D}x$, these sub-matrices
are column vectors of length $m$,
\[
P^{(1)} = p^{(1)}_1, \;\;
P^{(8)} = p^{(8)}_1, \;\; \mbox{and} \;\;
P^{(15)} = p^{(15)}_1 .
\]
The remaining sub-matrices are blurred versions of the left-integer and right-integer codes $L$ and $R$, represented as $m$-by-10 nonnegative real matrices.  We
write each of them as
\begin{equation}
P^{(j)} = \left[ p_1^{(j)} ~p_2^{(j)} ~\dots ~ p_{10}^{(j)} \right], \;\; j \ne 1, 8, 15,
\label{Pcolumns}
\end{equation}
where each $p_k^{(j)}$, $k =1,2, ..., 10$, is a column vector of 
length $m$. 

Recall that the over-sampling rate $r=m/n$ indicates the
number of time samples per minimal bar code width.
Given $r$, we can partition the rows of ${\cal P}$ into $15$ blocks, each block with index set $I_j$ of size 
  $|I_j|$, so that each sub-matrix is well-localized within a single block.  We know that if $P^{(1)}$ and $P^{(15)}$ correspond 
 to samples of the 3-bar sequence ``101'' or ``black-white-black'', so $|I_1| = |I_{15}| =
3r$.  The sub-matrix $P^{(8)}$ corresponds to samples from the middle 5 bar-sequence 
so $|I_8|=5r$.  Each remaining sub-matrix corresponds to samples from
a digit of length 7 bars, therefore $|I_j|=7r$ for $j\ne 1, 8, 15$.  

We can now give a quantitative measure describing how
`block-diagonal' the forward map is.    To this end, let $\varepsilon$
be the infimum of all $\epsilon > 0$ satisfying both
\begin{equation}\label{energyconc1}
\left\| \left. p_{k}^{(j)} \right|_{[m] \setminus I_j} \right\|_1 < \epsilon, 
\;\;\mbox{for all}\;\;
j \in [15], \; k \in [10],
\end{equation}
and
\begin{equation}\label{energyconc2}
\left\| \left. \left( \sum^{15}_{j' = j+1} p_{k_{j'}}^{(j')}
  \right) \right|_{I_j} \right \|_1 < \epsilon, \;\; \mbox{for
    all}\;\;
j \in [15], \; \mbox{and all choices of} \;\; k_{j+1}, \dots, k_{15} \in [10].
\end{equation}
The magnitude of $\varepsilon$ indicates to what extent the energy of each column 
of ${\cal P}$ is localized within its proper block.  If there were no blur, there would be 
no overlap between blocks and $\epsilon = 0$.  

\begin{figure}[h]
\centerline{
\includegraphics[width=5cm]{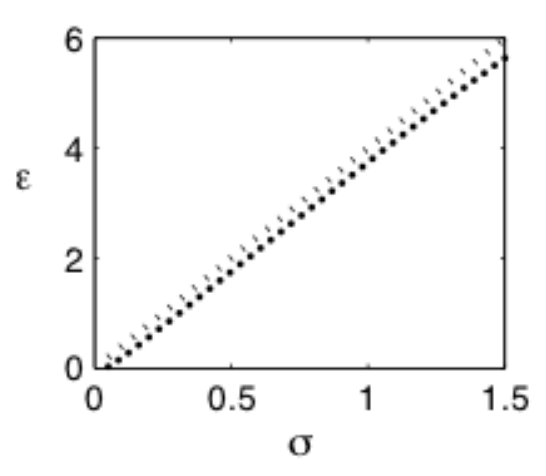} \hspace{1cm}
\includegraphics[width=5cm]{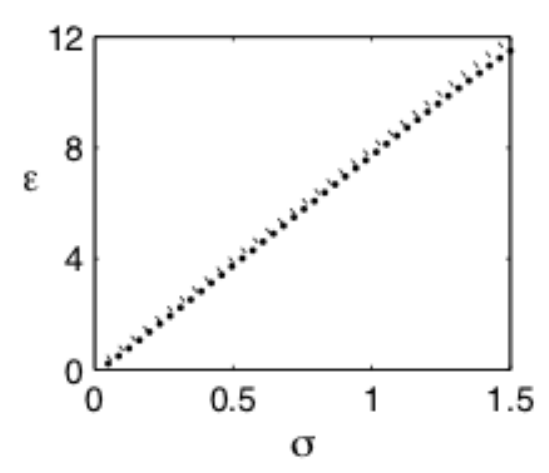}
}
\caption{For oversampling ratios $r=10$ (left) and $r=20$ (right) and $\alpha = 1$, the thick
line represents the minimal value of $\varepsilon$ 
  satisfying \eqref{energyconc1} and \eqref{energyconc2} in terms of $\sigma$.
   The thin line in each plot represents the function
  $(2/5)\sigma r$.}
\label{Fig:4}
\end{figure}

Simulation results such as those in Figure \ref{Fig:4} suggest that for $\alpha = 1$, the value of $\varepsilon$ in the forward map ${\cal P}$ can be expressed in terms of the oversampling ratio $r$ and Gaussian standard deviation $\sigma$ according to the formula $\varepsilon = (2/5) \sigma r$, at least over the relevant range of blur $0 \leq \sigma \leq 1.5$.  By linearity of the forward map with respect to the amplitude $\alpha$, this implies that more generally $\varepsilon = (2/5) \alpha \sigma r$.

\subsection{Column incoherence}
\begin{figure}[h]
\centerline{
\includegraphics[width=5cm]{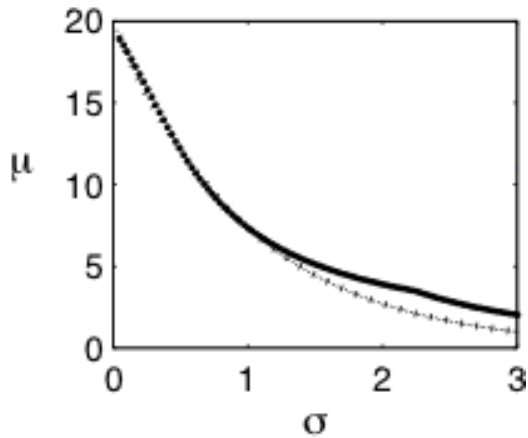}\hspace{1cm}
\includegraphics[width=5cm]{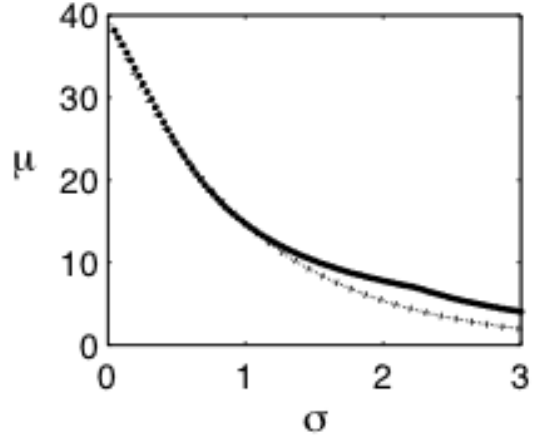}
}
\caption{For oversampling ratios $r=10$ (left) and $r=20$
  (right), we plot the minimal column separation $\mu = \min_{k_1
    \neq k_2} \left\| p_{k_1}^{(j)} - p_{k_2}^{(j)} \right\|_1$ for the 
    forward map ${\cal P} = {\cal G}(\sigma){\cal D}$, as a function of
     the standard deviation $\sigma$ of the Gaussian kernel.  The plots suggest that $\mu \approx 2 \alpha r e^{-\sigma} $ for $\sigma \leq 1$. }
\label{Fig:5}
\end{figure}

We now highlight another property of the forward map ${\cal P}$ that allows for robust
bar code reconstruction.  
The left-integer and right-integer codes for the UPC bar code, as enumerated
in Table \eqref{LRcodes}, are well-separated by design: the $\ell_1$-distance
 between any two distinct codes is greater than or equal to $2$.
 Consequently, if 
$D_k$  are the columns of the bar code dictionary
${\cal D}$, then $\min_{k_1 \neq k_2} \|D_{k_1} - D_{k_2} \|_1 = 2$.
This implies for the forward map ${\cal P} = \alpha {\cal G}(\sigma) {\cal D}$ that 
when there is no blur, i.e. $\sigma = 0$,
\begin{equation}
\label{min_separation}
\mu := \min_{j,k_1 \neq k_2} \left\| p_{k_1}^{(j)} - p_{k_2}^{(j)} \right\|_1 =
\min_{j,k_1 \neq k_2} \left\| \left. p_{k_1}^{(j)} \right|_{I_j} -
\left. p_{k_2}^{(j)}\right|_{I_j} \right\|_1 = 2 \alpha  r,
\end{equation}
where $r$ is the over-sampling ratio.
As the blur increases from zero, the column separation factor $\mu = \mu(\sigma, \alpha, r)$ decreases
smoothly.  In Figure \ref{Fig:5} we plot $\mu$ versus $\sigma$ for 
different oversampling ratios, as obtained from numerical simulation.  Simulations such as these suggest that $\mu$ closely follows the curve 
$\mu \approx 2 \alpha r e^{-\sigma}$, at least in the range $\sigma \leq 1$.  

\section{A simple decoding procedure for UPC bar codes}

We know from the bar code determination problem (\ref{noisydata}) that 
without additive noise, the observed data $d$ is the sum
of 15 columns from ${\cal P}$, one column from each block $P^{(j)}$.  Based 
on this observation, we will  employ a reconstruction algorithm which, once
initialized, selects the column from the successive block to minimize the norm of the data 
remaining after the column is subtracted.  This greedy
 algorithm is described in pseudo-code as follows.  

\begin{center}
\parbox{5in}{
{\bf Algorithm 1: Recover UPC Bar Code}\\
\rule[0pt]{5in}{0.5pt}

\tt
initialize: \\
\null\hspace{0.7cm} for $\ell=1, 62, 123$, $x_\ell=1$\\
\null\hspace{0.7cm} else $x_\ell=0$\\
\null\hspace{0.7cm} $\delta \leftarrow d$\\
for $j=2:7,~9:14$\\
\null\hspace{.7cm} $k_{\min} = \argmin_k  \left\| \delta - p_k^{(j)} \right\|_1$\\
\null\hspace{.7cm} if $j\leq 7$, $\ell\leftarrow 1+10(j-2)+k_{\min}$\\
\null\hspace{.7cm} else $\ell\leftarrow 62+10(j-9)+k_{\min}$\\
\null\hspace{.7cm} $x_\ell \leftarrow 1$\\
\null\hspace{.7cm} $r \leftarrow \delta - p^{(j)}_{k_{\min}}$\\
end\\
\rule[0pt]{5in}{0.5pt}
}
\end{center}

%

\section{Analysis of the algorithm}
\label{sec:AlgAnal}

Algorithm 1 recovers the bar code one digit at a time by iteratively scanning through the observed data.  The runtime complexity of the method is dominated by the $12$ calculations of $k_{\rm min}$ performed by the algorithm's single loop over the course of its execution.  Each one of these calculations of $k_{\rm min}$ consists of $10$ computations of the $\ell_1$-norm of a vector of length $m$.  Thus, the runtime complexity of the algorithm is $O(m)$, and can be executed in less than a second for standard UPC bar code proportions.\footnote{In practice, when $\sigma$ is not too large, a `windowed' vector of length less than $m$ can be used to approximate $\left\| \delta - p_k^{(j)} \right\|_1$ for each $k,j$.  This can reduce the constant of proportionality associated with the runtime complexity.}  

\subsection{Recovery of the unknown bar code}
Recall that
the 12 unknown digits in the unknown bar code $c$ are represented by the sparse vector $x$
 in $c = {\cal D}x$. We already know that $x_1=x_{62}=x_{123} =1$ as these
elements corresponds to the mandatory start, middle, and end sequences.
Assuming for the moment that the forward map 
${\cal P}$ is \emph{known}, i.e., that both $\sigma$ and $\alpha$ are known, we now prove that the greedy algorithm
will reconstruct the correct bar code from noisy data $d = {\cal P}x + h$ as long as ${\cal P}$ is sufficiently block-diagonal and if its columns are sufficiently incoherent.   In
the next section we will extend the analysis to the case where $\sigma$ and $\alpha$
are unknown.

\newtheorem{Theorem}{Theorem}
\begin{Theorem}
Suppose $I_1, \dots, I_{15} \subset [m]$ and 
$\varepsilon \in \mathbb{R}$ satisfy the conditions 
(\ref{energyconc1})-(\ref{energyconc2}).  Then, Algorithm~1 
will correctly recover a bar code signal $x$ from noisy data $d = {\cal P}x + h$ provided that
\begin{equation}
\left\| \left. p_{k_1}^{(j)} \right|_{I_j} - \left. p_{k_2}^{(j)} \right|_{I_j} \right\|_1 > 
2 \left( \left\| h |_{I_j} \right\|_1 + 2\varepsilon \right)
\label{VictoryCondition}
\end{equation} 
for all $j \in [15]$ and $k_1, k_2 \in [10]$ with $k_1 \neq k_2$.
\end{Theorem}

\noindent \textit{Proof:}\\

Suppose that 
$$d = {\cal P} x + h = \sum^{15}_{j=1} p_{k_j}^{(j)} + h.$$
Furthermore, denoting $k_j = k_{\min}$ in the for-loop in Algorithm 1,
suppose that $k_2, \dots, k_{j'-1}$ have already been correctly recovered.
Then the 
residual data, $\delta$, at this stage of the algorithm will be 
$$\delta = p_{k_{j'}}^{(j')} + \delta_{j'} + h,$$
where $\delta_{j'}$ is defined to be 
$$\delta_{j'} = \sum^{15}_{j = j' + 1} p_{k_j}^{(j)}.$$
We will now show that the ${j'}$th execution of the for-loop
will correctly recover $p_{k_{j'}}^{(j')}$, thereby establishing the desired result by
induction.

Suppose that the ${j'}$th execution of the for-loop incorrectly recovers $k_{\rm err} 
\neq k_{j'}$.  This happens if
$$
\left\| \delta - p_{k_{\rm err}}^{(j')} \right\|_1 \leq \left\| \delta - p^{(j')}_{k_{j'}} \right\|_1.
$$
In other words, we have that
\begin{align}
\left\| \delta - p_{k_{\rm err}}^{(j')} \right\|_1 & = 
  \left\| \delta \big|_{I_{j'}}  - \left. p_{k_{\rm err}}^{(j')} \right|_{I_{j'}} \right\|_1 
+ \left\| \delta \big|_{I^{\rm c}_{j'}} - \left. p_{k_{\rm err}}^{(j')} \right|_{I^{\rm c}_{j'}} \right\|_1 \nonumber \\ 
& \geq
\left\| \left. p^{(j')}_{k_{j'}} \right|_{I_{j'}} - \left. p_{k_{\rm err}}^{(j')} \right|_{I_{j'}} \right\|_1 
- \left\| \delta_{j'} \big|_{I_{j'}} \right\|_1 - \left\| h \big|_{I_{j'}} \right\|_1 \nonumber \\
&+ \left\| \delta_{j'} \big|_{I^{\rm c}_{j'}}  + h \big|_{I^{\rm c}_{j'}} \right\|_1 
 -  \left\| \left. p^{(j')}_{k_j'} \right|_{I^{\rm c}_{j'}} \right\|_1 
 -  \left\| \left. p_{k_{\rm err}}^{(j')} \right|_{I^{\rm c}_{j'}} \right\|_1  \nonumber \\ 
& \geq
\left\| \left. p_{k_{j'}}^{(j')} \right|_{I_{j'}} - \left. p_{k_{\rm err}}^{(j')} \right|_{I_{j'}} \right\|_1 
+ \left\|  \delta_{j'} \big|_{I^{\rm c}_{j'}} + h \big|_{I^{\rm c}_{j'}} \right\|_1 
- \left\| h \big|_{I_{j'}} \right\|_1 - 3\varepsilon \nonumber
\end{align}
from conditions (\ref{energyconc1}) and (\ref{energyconc2}).  
To finish, we simply simultaneously add and
subtract $\| \delta_{j'} |_{I_{j'}} + h|_{I_{j'}} \|_1$ 
from the last expression to arrive at a contradiction to the supposition that 
$k_{\rm err} \neq k_{j'}$:
\begin{align}
\left\| \delta - p_{k_{\rm err}}^{(j')} \right\|_1 &\geq 
\left( \left\| \left. p_{k_{j'}}^{(j')} \right|_{I_{j'}} - 
   \left. p_{k_{\rm err}}^{(j')} \right|_{I_{j'}} \right\|_1 
- 2 \left\| h \big|_{I_{j'}} \right\|_1 - 4 \varepsilon \right) 
+ \left\| \delta_{j'} + h \right\|_1 \nonumber \\
&= \left( \left\| \left. p_{k_{j'}}^{(j')} \right|_{I_{j'}} - 
\left. p_{k_{\rm err}}^{(j')} \right|_{I_{j'}} \right\|_1 - 
2 \left\| h \big|_{I_{j'}} \right\|_1 - 4 \varepsilon \right) 
+ \left\| \delta - p^{(j')}_{k_{j'}} \right\|_1 \nonumber \\
&> \left\| \delta - p^{(j')}_{k_{j'}} \right\|_1.
\label{LowerBound}
\end{align}
\hfill{$\Box$}

\begin{rek}
\label{rek:empirical}
\textnormal{
Equation~\eqref{energyconc1} implies that 
$$\min_{j, k_1 \neq k_2} \left\| \left. p_{k_1}^{(j)} \right|_{I_j} 
  - \left. p_{k_2}^{(j)} \right|_{I_j} \right\|_1 ~\geq~ \min_{j, k_1 \neq k_2} \left\| p_{k_1}^{(j)} - p_{k_2}^{(j)}
\right\|_1 - 2\varepsilon ~=~ \mu - 2\varepsilon.\footnote{See equation~\eqref{min_separation} for the definition of $\mu$.}$$
Thus, the recovery condition \eqref{VictoryCondition} in Theorem 1 will hold whenever
\[
\mu - 2\varepsilon > 2 \left( \left\| h |_{I_j} \right\|_1 + 2\varepsilon \right).
\]
Using the empirical relationships $\varepsilon = (2/5) \alpha r \sigma$ and $\mu = 2\alpha r e^{-\sigma}$, we obtain the following
upper bound on the level of sufficient noise for successful recovery:
\begin{equation}
\label{VictoryCond2}
\max_{j \in [12]} \left\| h|_{I_j} \right\|_1 < \alpha r (e^{-\sigma} - (6/5) \sigma).
\end{equation}
In practice the Gaussian blur width $2\sqrt{2\ln(2)}\sigma$ does not exceed 
the minimum width of the bar code, which we have normalized 
to be $1$.  This translates to a maximal standard deviation of  $\sigma
\approx .425$, and a noise ceiling in \eqref{VictoryCond2} of 
\begin{equation}
\label{noisetol}
\max_{j \in [12]} \left\| h |_{I_j} \right\|_1 \leq .144 \alpha r.
\end{equation}
This should be compared to the $\ell_1$-norm of the bar code signal over a block; the average $\ell_1$ norm between the left-integer and right-integer codes is $3.5 \alpha$.
}
\end{rek}

\begin{rek}
\textnormal{
In practice it may be beneficial to apply Algorithm 1 several times, 
each time changing the order in which the
digits are decoded.  For example, if the distribution of the noise is known
in advance, it would be beneficial to to initialize the algorithm 
 in regions of the bar code with less noise.}
\end{rek}

\subsection{Stability of the greedy algorithm with respect to parameter estimation}

\subsubsection*{Insensitivity to unknown $\alpha$}

In the previous section we assumed a known Gaussian convolution
matrix $\alpha {\cal G}(\sigma)$.  In fact, this is
generally not the case.  In practice both
$\sigma$ and $\alpha$ must be estimated since these parameters depend
on the distance from the
scanner to the bar code, the reflectivity of the scanned
surface, the ambient light, etc.  This
means that in practice, Algorithm 1 will be decoding bar codes using only an
approximation to $\alpha {\cal G}(\sigma)$.
Suppose that the true standard deviation generating a sampled sequence $d$ is $\sigma$, but that Algorithm 1 uses a different value $\widehat{\sigma}$ for reconstruction.  
We can regard the error incurred by $\widehat{\sigma}$ as additional additive noise in our sensitivity analysis, setting $h' =  h + \alpha \big({\cal G}(\sigma) - {\cal G}(\widehat{\sigma}) \big) {\cal D} x$ and rewriting the inverse problem as
\begin{eqnarray} \label{noisydata0}
d &=& \alpha {\cal G}(\sigma){\cal D}x + h \nonumber \\
 &=& \alpha {\cal G}(\widehat{\sigma}){\cal D}x + \Big(h + \alpha \big({\cal G}(\sigma) - {\cal G}(\widehat{\sigma}) \big) {\cal D} x \Big) \nonumber \\
 &=&  \alpha {\cal G}(\widehat{\sigma}){\cal D}x + h'.
\end{eqnarray}
We now describe a procedure for estimating $\alpha$.  Note that the middle portion of the observed data of length $5r$, $\left. d_{mid} = d \right|_{I_8}$, represents a blurry image of the known middle pattern $M = [01010]$.  Let ${\cal P} = {\cal G}(\widehat{\sigma}){\cal D}$ be the forward map generated by the estimate $\widehat{\sigma}$ when $\alpha = 1$, and consider the sub-matrix
$$p_{mid} = \left. P^{(8)} \right|_{I_8}$$
which is also a vector of length $5r$.  If $\widehat{\sigma} = \sigma$ or $\widehat{\sigma} \approx \sigma$,\footnote{Here we have assumed that the noise level is low.  In noisier settings it should be possible to develop more effective methods for estimating $\alpha$ by incorporating the characteristics of the scanning noise.} we expect a good estimate for $\alpha$ to be the least squares solution
\begin{equation}
\label{alphals}
\widehat{\alpha}= \arg \min_{a} \| a p_{mid} - d_{mid} \|_2 = p_{mid}^T d_{mid}/\|\ p_{mid} \|_2^2.
\end{equation}

Dividing both sides of the equation \eqref{noisydata0} by $\widehat{\alpha}$, the inverse problem becomes
\begin{eqnarray} \label{inv:a}
\frac{d}{\widehat{\alpha}}
&=& \frac{\alpha}{\widehat{\alpha}} {\cal G}(\widehat{\sigma}){\cal D} x +
\frac{1}{\widehat{\alpha}} h'.
\end{eqnarray}

Suppose that
$1 - \gamma \leq \alpha/\widehat{\alpha} \leq 1 + \gamma$ for some
$0<\gamma<1$.  Then fixing the data to be $\widehat{d} = d/\widehat{\alpha}$ and
fixing forward map
to be ${\cal P} = {\cal G}(\widehat{\sigma}){\cal D}$,
the recovery conditions \eqref{energyconc1}, \eqref{energyconc2}, and \eqref{VictoryCondition} become respectively
\begin{enumerate}{\itemsep=0cm\parsep=0cm}
\item $ \left\| \left. p_{k}^{(j)} \right|_{[m] \setminus
  I_j} \right\|_1 < \frac{\varepsilon}{1+\gamma}$ for all $j \in [15]$ and $k \in [10]$. 
\item $ \left\| \left. \left( \sum^{15}_{j' = j+1} 
p_{k_{j'}}^{(j')} \right)\right|_{I_j}
\right \|_1 < \frac{\varepsilon}{1+\gamma}$ for all $j \in [14]$ and valid 
$k_{j'} \in [10]$.
\item $\left\| \left. p_{k_1}^{(j)} \right|_{I_j} 
- \left. p_{k_2}^{(j)} \right|_{I_j} \right\|_1 > 
2 \Big( \frac{1}{\alpha} \left\| \left. h \right|_{I_j} \right\|_1 + \left\| \left. \big( {\cal G}(\sigma) - {\cal G}(\widehat{\sigma})  \big) {\cal D} x  \right|_{I_j} \right\|_1 + \frac{2\varepsilon}{1-\gamma} \Big) $
\end{enumerate}
Consequently, if $\sigma \approx \widehat{\sigma}$ and $1 \lessapprox \alpha \approx \widehat{\alpha}$,  the conditions for correct bar code reconstruction do not change much.

\subsubsection*{Insensitivity to unknown $\sigma$}

We have seen that one way to estimate the scaling
$\alpha$ is to guess a value for $\sigma$ and perform a least-squares
fit of the observed data.  In doing so, we found that the sensitivity of the 
recovery process with respect to $\sigma$ is proportional to the quantity
\begin{equation}
\label{Gterm}
\left\| 
\left. \left( {\cal G}(\sigma) - {\cal G}(\widehat{\sigma}) \right)
{\cal D} x \right|_{I_j} \right\|_1
\end{equation}
in the third condition immediately above.  Note that all the entries of the matrix ${\cal G}(\sigma) -
{\cal G}(\widehat{\sigma})$ will be small whenever $\widehat{\sigma}
\approx \sigma$.  Thus, Algorithm~1 should be able to tolerate small
parameter estimation errors as long as the ``almost'' block diagonal
matrix formed using $\widehat{\sigma}$ exhibits a sizable difference
between any two of its digit columns which might (approximately)
appear in any position of a given UPC bar code.

To get a sense of the size of the term \eqref{Gterm}, let us further investigate the expressions involved.
Recall that using the dictionary matrix ${\cal D}$, a bar code sequence of 0's and 1's is
given by $c={\cal D}x$.  When put together with the bar code function representation
(\ref{barcoderep}), we see that
\[
\left[ {\cal G}(\sigma) {\cal D} x \right]_i =
\int g_\sigma (t_i -t) z(t) dt,
\]
where
\[
g_\sigma (t) = \frac{1}{\sqrt{2\pi}\sigma} e^{-(\frac{t^2}{2\sigma^2})} .
\]
Therefore, we have
\begin{equation}
\label{GDx}
\left[ {\cal G}(\sigma) {\cal D} x \right]_i =
\sum_{j=1}^n c_j \int_{j-1}^j g_\sigma (t_i-t) dt .
\end{equation}
Now, using the definition for the cumulative distribution function for normal
distributions
\[
\Phi(x)=\frac{1}{\sqrt{2\pi}} \int_{-\infty}^x e^{-t^2/2} dt,
\]
we see that
\[
\int_{j-1}^j g_\sigma (t_i-t) dt = \Phi\left( \frac{t_i-j+1}{\sigma} \right) -
                                \Phi\left( \frac{t_i-j}{\sigma} \right) .
\]
and we can now rewrite (\ref{GDx}) as
\[
\left[ {\cal G}(\sigma) {\cal D} x \right]_i =
\sum_{j=1}^n c_j \left[ \Phi\left( \frac{t_i-j+1}{\sigma} \right) -
                \Phi\left( \frac{t_i-j}{\sigma} \right) \right] .
\]
We now isolate the term we wish to analyze:
\begin{align*}
& \left[ \left( {\cal G}(\sigma) - {\cal G}(\widehat{\sigma}) \right) {\cal D} x \right]_i\\
& \hspace{0.5cm} = \sum_{j=1}^n c_j \left[ \Phi\left( \frac{t_i-j+1}{\sigma} \right) -
\Phi\left( \frac{t_i-j+1}{\widehat{\sigma}} \right) -                                 
\Phi\left( \frac{t_i-j}{\sigma} \right) +
\Phi\left( \frac{t_i-j}{\widehat{\sigma}} \right) \right] .
\end{align*}
We are interested in the error
\begin{eqnarray}
&& \left|
 \left[ \left( {\cal G}(\sigma) - {\cal G}(\widehat{\sigma}) \right) {\cal D} x \right]_i \right| \nonumber \\
&\leq& \sum_{j=1}^n c_j \left| \Phi\left( \frac{t_i-j+1}{\sigma} \right) -
\Phi\left( \frac{t_i-j+1}{\widehat{\sigma}} \right) -                                 
\Phi\left( \frac{t_i-j}{\sigma} \right) +
\Phi\left( \frac{t_i-j}{\widehat{\sigma}} \right) \right|  \nonumber \\
&\leq& \sum_{j=1}^n \left| \Phi\left( \frac{t_i-j+1}{\sigma} \right) -
\Phi\left( \frac{t_i-j+1}{\widehat{\sigma}} \right) \right| +
                                 \left| \Phi\left( \frac{t_i-j}{\sigma} \right) -
\Phi\left( \frac{t_i-j}{\widehat{\sigma}} \right) \right| , \nonumber \\
&\leq& 2\sum_{j=0}^{n} \left| \Phi\left( \frac{t_i-j}{\sigma} \right) -
\Phi\left( \frac{t_i-j}{\widehat{\sigma}} \right) \right|. \nonumber
\end{eqnarray}

Suppose that $\xi = (\xi_k)$ is the vector of values $| t_i - j|$ for fixed $i$, running $j$, sorted in order of increasing magnitude.  Note that $\xi_1$ and $\xi_2$ are less than or equal to 1, and  $\xi_3 \leq \xi_1 + 1$, $\xi_4 \leq \xi_2 + 1$, and so on.  We can center the previous bound around $\xi_1$ and $\xi_2$, giving 
\begin{eqnarray}
\left|
 \left[ \left( {\cal G}(\sigma) - {\cal G}(\widehat{\sigma}) \right) {\cal D} x \right]_i \right| 
\leq \sum_{j=0}^{n} \left| \Phi\left( \frac{\xi_1 + j}{\sigma} \right) -
\Phi\left( \frac{\xi_1 + j}{\widehat{\sigma}} \right) \right| +  \left| \Phi\left( \frac{\xi_2 + j}{\sigma} \right) -
\Phi\left( \frac{\xi_2 + j}{\widehat{\sigma}} \right) \right|.
\label{eq:terms}
\end{eqnarray}

Next we simply majorize the expression
\[
\label{tomajorize}
f(x) = \Phi\left( \frac{x}{\sigma} \right) - \Phi \left( \frac{x}{\widehat{\sigma}} \right).
\]
To do so, we take the derivative and find the critical points, which turn out to be
\[
x_* = \pm \sqrt{2}\sigma\widehat{\sigma} 
\sqrt{ \frac{\log \sigma - \log \widehat{\sigma}}{\sigma^2-\widehat{\sigma}^2} }.
\]
Therefore, each term in the summand \eqref{eq:terms} can be bounded by
\begin{eqnarray}
  \left| \Phi\left( \frac{\xi + j}{\sigma} \right) -
\Phi\left( \frac{\xi + j}{\widehat{\sigma}} \right) \right|
&\leq&   \left| \Phi\left( \sqrt{2}\widehat{\sigma} 
\sqrt{ \frac{\log \sigma - \log \widehat{\sigma}}{\sigma^2-\widehat{\sigma}^2} } \right) -
\Phi\left( \sqrt{2}\sigma 
\sqrt{ \frac{\log \sigma - \log \widehat{\sigma}}{\sigma^2-\widehat{\sigma}^2} } \right) \right| \nonumber \\
&:=& \triangle_1(\sigma, \widehat{\sigma}) . \label{bound1}
\end{eqnarray}

On the other hand, the terms in the sum decrease exponentially as $j$ increases. 
To see this, recall the simple bound
\begin{align*}
1 - \Phi(x) \hspace{2mm} =& \hspace{2mm} \frac{1}{\sqrt{2\pi}} \int_{x}^{\infty} e^{-t^2/2} dt \hspace{2mm} \leq \hspace{2mm}   \frac{1}{\sqrt{2\pi}} \int_{x}^{\infty} \frac{t}{x} e^{-t^2/2} dt \hspace{2mm} = \hspace{2mm} \frac{e^{-x^2/2}}{x\sqrt{2\pi}}.
\end{align*}
Writing $\sigma_{max} = \max\{(\sigma, \widehat{\sigma})\}$, and noting that $\Phi(x)$ is a positive, increasing function, we have for $\xi \in [0,1)$
\begin{eqnarray}
  \left| \Phi\left( \frac{\xi + j}{\sigma} \right) -
\Phi\left( \frac{\xi + j}{\widehat{\sigma}} \right) \right|
&\leq&  1 - \Phi\left( \frac{\xi + j}{\sigma_{max} } \right)  \nonumber \\
 &\leq& \frac{\sigma_{max}}{(\xi + j) \sqrt{2\pi}}e^{-(\xi + j)^2/(2\sigma_{max}^2)}  \nonumber \\
&\leq& \frac{\sigma_{max}}{ (\xi + j) \sqrt{2\pi}}e^{-(\xi + j)/(2\sigma_{max}^2)} \quad \text{ if } j \geq 1 \nonumber \\
&=& \frac{\sigma_{max}}{(\xi + j) \sqrt{2\pi}}\Big( e^{-(2\sigma_{max}^2)^{-1}}\Big)^{\xi + j }  \nonumber \\
&\leq& \frac{\sigma_{max}}{j \sqrt{2\pi}}\Big( e^{-(2\sigma_{max}^2)^{-1}}\Big)^{j} \nonumber  \\
&:=& \Delta_2(\sigma_{max},j). \label{bound2}
\end{eqnarray} 
Combining the bounds \eqref{bound1} and \eqref{bound2}, 
$$
  \left| \Phi\left( \frac{\xi + j}{\sigma} \right) -
\Phi\left( \frac{\xi + j}{\widehat{\sigma}} \right) \right| \leq \min{\Big( (\triangle_1(\sigma, \widehat{\sigma}), \triangle_2(\sigma_{max}, j) \Big)}.
$$
Suppose that $j_1$ is the smallest integer in absolute value such that $\triangle_2(\sigma_{max}, j_1) \leq \triangle_1(\sigma, \widehat{\sigma})$. Then from this term on, the summands in \eqref{eq:terms} can be bounded by a geometric series:
\begin{eqnarray}
 \sum_{ j \geq j_1}^n   \left| \Phi\left( \frac{\xi + j}{\sigma} \right) -
\Phi\left( \frac{\xi + j}{\widehat{\sigma}} \right) \right| 
&\leq&  \frac{2\sigma_{max}}{ j_1\sqrt{2\pi}} \sum_{j \geq j_1}  a^j, \quad a =  e^{-(2\sigma_{max}^2)^{-1}}  \nonumber \\
&\leq&  \frac{2\sigma_{max}}{ j_1 \sqrt{2\pi}} \cdot a^{j_1} (1 - a)^{-1} . \nonumber
\end{eqnarray}
We then arrive at the bound
\begin{eqnarray}
\label{boundG1}
 \left|
 \left[ \left( {\cal G}(\sigma) - {\cal G}(\widehat{\sigma}) \right) {\cal D} x \right]_i \right| &\leq& 2\cdot  j_1 \triangle_1(\sigma, \widehat{\sigma}) +  \frac{4\sigma_{max}\cdot a^{j_1}(1-a)^{-1}}{j_1\sqrt{2\pi}} \nonumber \\
 &=:& B(\sigma, \widehat{\sigma}).
\end{eqnarray}
The term \eqref{Gterm} can then be bounded according to
\begin{eqnarray}
\label{boundG}
\left\| 
\left. \left( {\cal G}(\sigma) - {\cal G}(\widehat{\sigma}) \right)
{\cal D}x \right|_{I_j} \right\|_1 &\leq& | I_j |  B(\sigma, \widehat{\sigma}) \leq 7r B(\sigma, \widehat{\sigma}),
\end{eqnarray}
where $r = m/n$ is the over-sampling rate.

Recall that in practice the width $2\sqrt{2\ln(2)}\sigma$ of the Gaussian kernel is
on the order of the minimum bar width in the bar code, which we normalized to $1$.  When the blur exactly equals the minimum bar width, we arrive at $\sigma
\approx .425$.  Below, we compute the error bound $B(\sigma,
\widehat{\sigma})$ for $\sigma = .425$ and several values of
$\widehat{\sigma}$.  

\medskip
\begin{center}
\begin{tabular}{|c|l|l|l|l|l|l|}\hline
$\widehat{\sigma}$ & 0.2 & 0.4 & .5 & 0.6 & .8  \\ \hline
$B(.425,\widehat{\sigma})$ & .3453 & .0651& .1608 &  .3071 &.589 \\ \hline
\end{tabular}
\end{center}

\medskip
While the bound \eqref{boundG} is very rough, note that the tabulated error bounds incurred by inaccurate $\sigma$ are at least roughly the same order of magnitude as the empirical noise level tolerance for the greedy algorithm, as discussed in Remark $6.1$.

\section{Numerical Evaluation}

In this section we illustrate with numerical examples the robustness
of the greedy algorithm to signal noise and imprecision in the
$\alpha$ and $\sigma$ parameter estimates.  We assume that neither
$\alpha$ nor $\sigma$ is known a priori, but that we have an estimate
$\widehat{\sigma}$ for $\sigma$.  We then compute an estimate
$\widehat{\alpha}$ from $\widehat{\sigma}$ by solving the
least-squares problem \eqref{alphals}.  

\begin{figure}[ht]

\subfigure[True parameter values: $\sigma = .45$, $\alpha = 1$.]{
\resizebox{0.48\columnwidth}{!}{\includegraphics{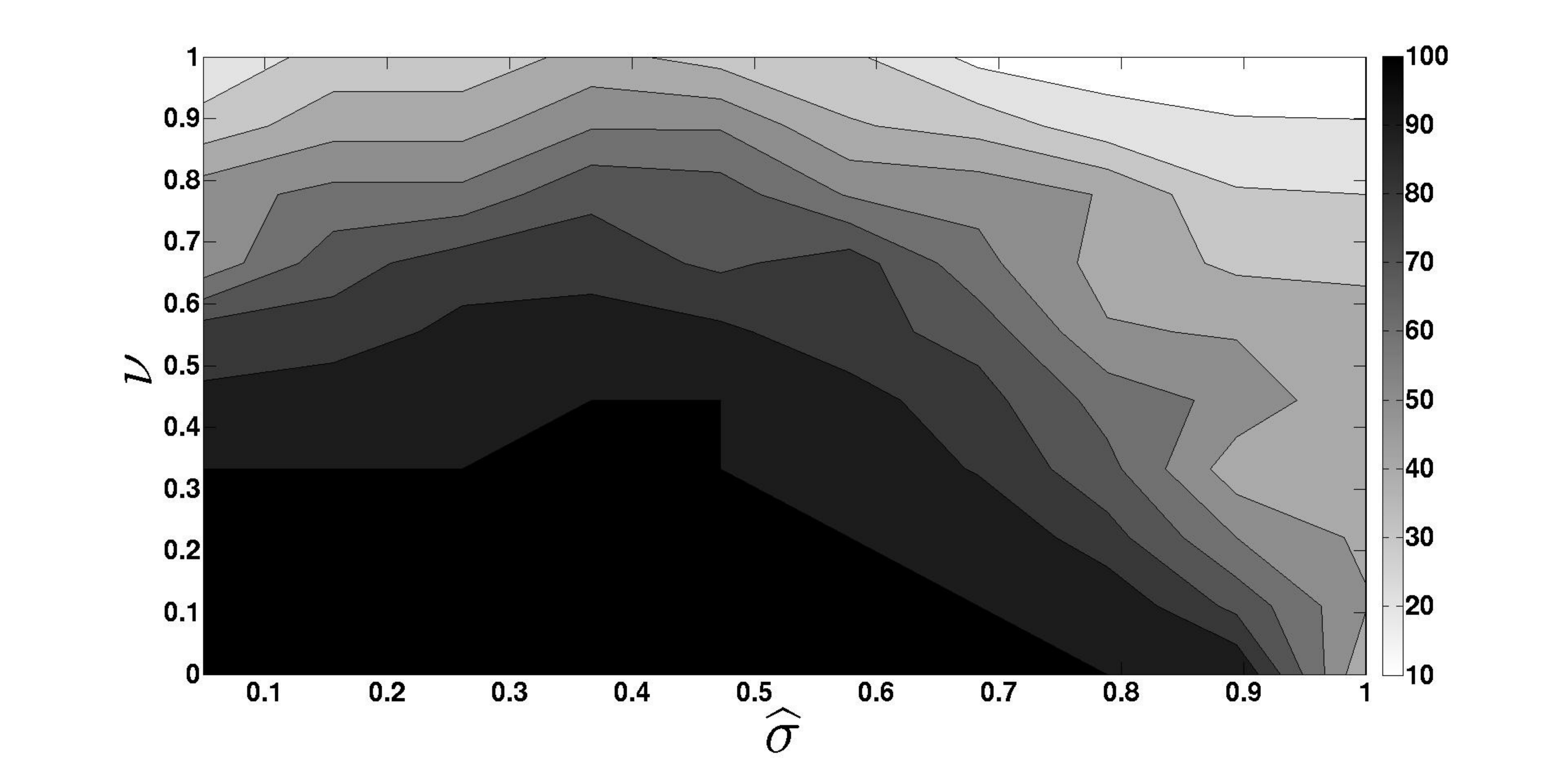}}
}
\subfigure[True parameter values: $\sigma = .75$, $\alpha = 1$]{
\resizebox{0.48\columnwidth}{!}{\includegraphics{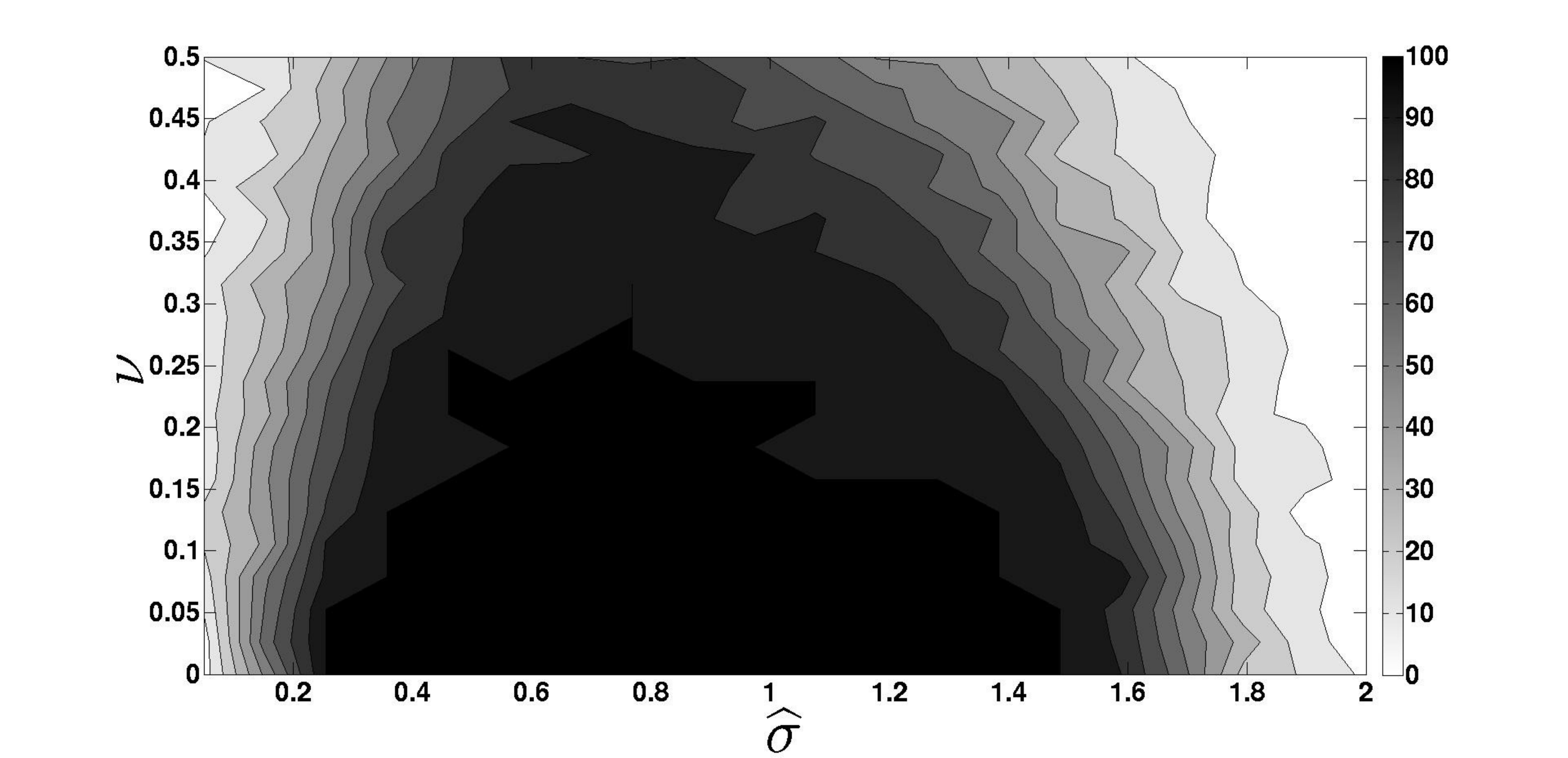}}
}
\caption{Recovery Probabilities when $\alpha = 1$ for two true
  $\sigma$ settings.  The shade in each phase diagram corresponds to
  the probability that the greedy algorithm will correctly recover a
  randomly selected bar code, as a function of the relative
  noise-to-signal level, $\nu = \frac{\| h \|_2}{\| \alpha {\cal
      G}(\sigma){\cal D} x \|_2}$, and the $\sigma$ estimate,
  $\widehat{\sigma}$.  Black represents correct bar code recovery with
  probability 1, while pure white represents recovery with probability
  0.  Each data point's shade (i.e., probability estimate) is based on
  100 random trials.}
\label{fig:relphase1}
\end{figure}

The phase diagrams in Figure~\ref{fig:relphase1} demonstrate the
insensitivity of the greedy algorithm to relatively large amounts of
noise.  These diagrams were constructed by executing Algorithm 1 on many trial input signals of
the form $d = \alpha {\cal G}(\sigma){\cal D} x + h$, where $h$ is mean zero
Gaussian noise. More specifically, each trial signal, $d$, was formed
as follows: a 12 digit number was first generated uniformly at
random, and its associated blurred bar code, $\alpha {\cal
  G}(\sigma){\cal D} x$, was formed using the oversampling ratio $r = m/n =
10$.  Second, a noise vector $n$ with independent and identically distributed entries $n_j \sim {\cal N}(0,1)$ was
generated and then rescaled to form the additive noise vector $h = \nu
\| \alpha {\cal G}(\sigma){\cal D} x \|_2 \left( n / \| n \|_2 \right)$.
Hence, the parameter $\nu = \frac{\| h \|_2}{\| \alpha {\cal
    G}(\sigma){\cal D} x \|_2}$ represents the noise-to-signal ratio of each
trial input signal $d$. 

We note that in laser-based scanners, there are two major sources of
noise: electronic noise \cite{kogan2008electronic}, which is often
modeled as $1/f$ noise \cite{dutta1981low}, and speckle noise
\cite{marom-kresic}, caused by the roughness of the paper. However,
the Gaussian noise used in our numerical experiments is sufficient for
the purpose of this work.

To create the phase diagrams in Figure~\ref{fig:relphase1}, the greedy
recovery algorithm was run on 100 independently generated trial input
signals for each of at least 100 equally spaced
$(\widehat{\sigma},\nu)$ grid points (a $10 \times 10$ mesh was used
for Figure~\ref{fig:relphase1}(a), and a $20 \times 20$ mesh for
Figure~\ref{fig:relphase1}(b)).  The number of times the greedy
algorithm successfully recovered the original UPC bar code determined
the color of each region in the $(\widehat{\sigma},\nu)$-plane.  The
black regions in the phase diagrams indicate regions of parameter
values where all $100$ of the $100$ randomly generated bar codes were
correctly reconstructed. The pure white parameter regions indicate
where the greedy recovery algorithm failed to correctly reconstruct
any of the $100$ randomly generated bar codes. 

Looking at Figure~\ref{fig:relphase1} we can see that the greedy
algorithm appears to be highly robust to additive noise.  For example,
when the $\sigma$ estimate is accurate (i.e., when $\widehat{\sigma}
\approx \sigma$) we can see that the algorithm can tolerate additive
noise with Euclidean norm as high as $0.25 \| \alpha {\cal
  G}(\sigma){\cal D}x \|_2$.  Furthermore, as $\widehat{\sigma}$
becomes less accurate the greedy algorithm's accuracy appears to
degrade smoothly.

\begin{figure}[ht]
\subfigure[True parameter values: $\sigma = .45$, $\alpha = 1$.]{
\resizebox{0.48\columnwidth}{!}{\includegraphics{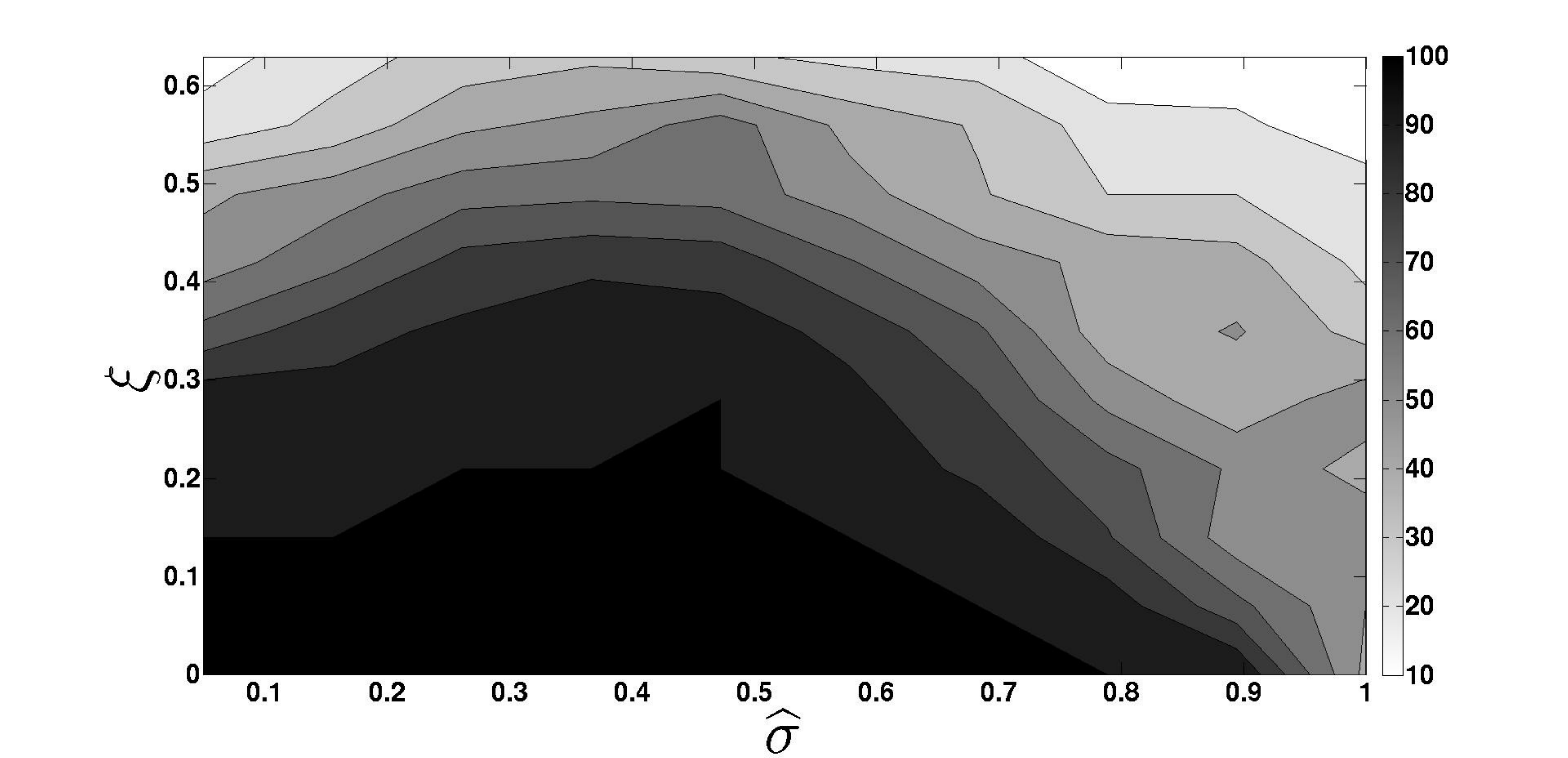}}
}
\subfigure[True parameter values: $\sigma = .75$, $\alpha = 1$]{
\resizebox{0.48\columnwidth}{!}{\includegraphics{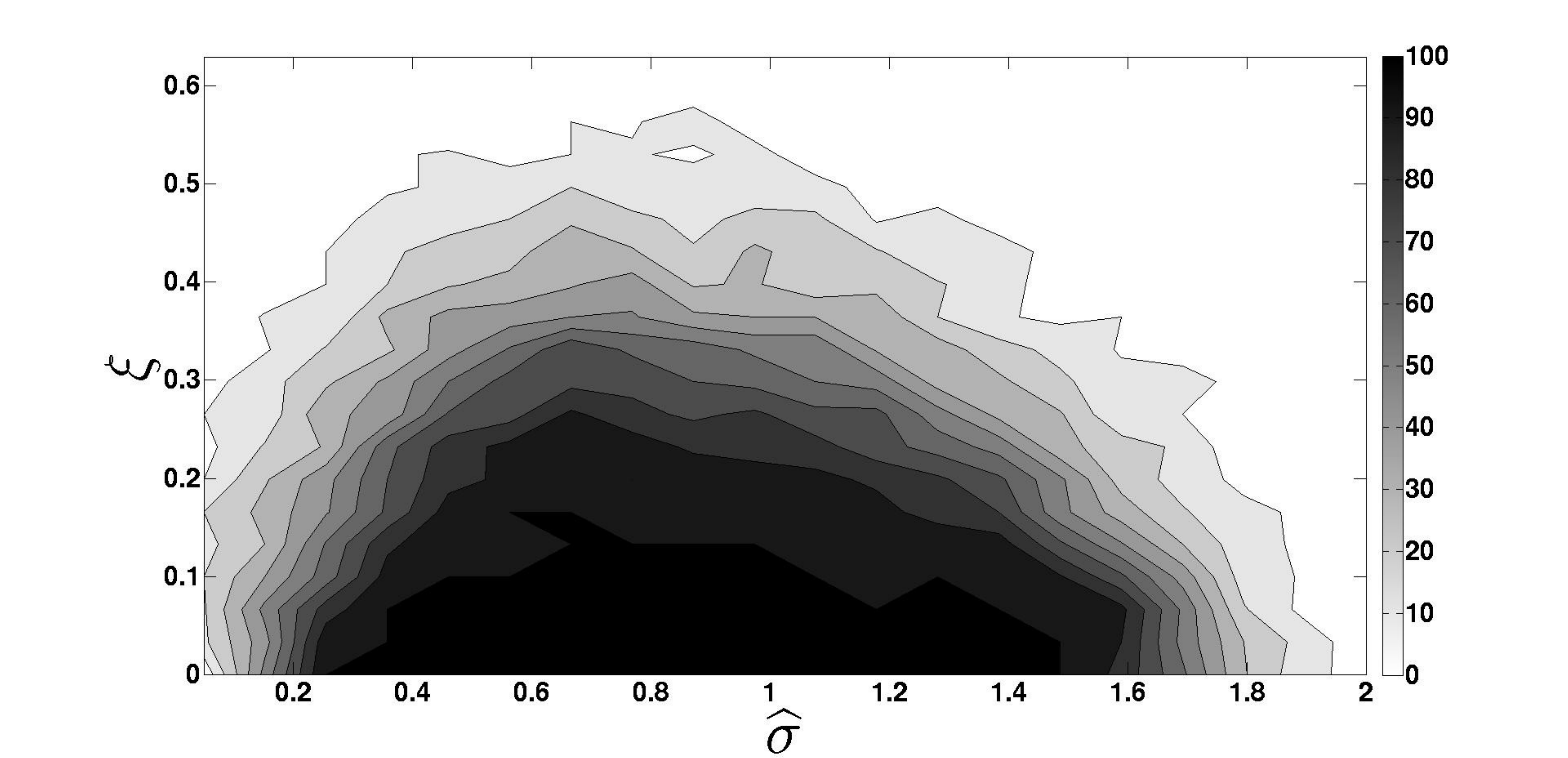}}
}
\caption{Recovery probabilities when $\alpha = 1$ for two true $\sigma$ settings.  The shade in each phase diagram corresponds to the probability that the greedy algorithm correctly recovers a randomly selected bar code, as a function of the additive noise standard deviation, $\xi$, and the $\sigma$ estimate, $\widehat{\sigma}$.  Black represents correct bar code recovery with probability 1, while pure white represents recovery with probability 0.  Each data point's shade (i.e., probability estimate) is based on 100 random trials.} 
\label{fig:phase1}
\end{figure}

\begin{figure}[ht]
\subfigure[True parameter values: $\sigma = .45$, $\alpha = .25$.]{
\resizebox{0.48\columnwidth}{!}{\includegraphics{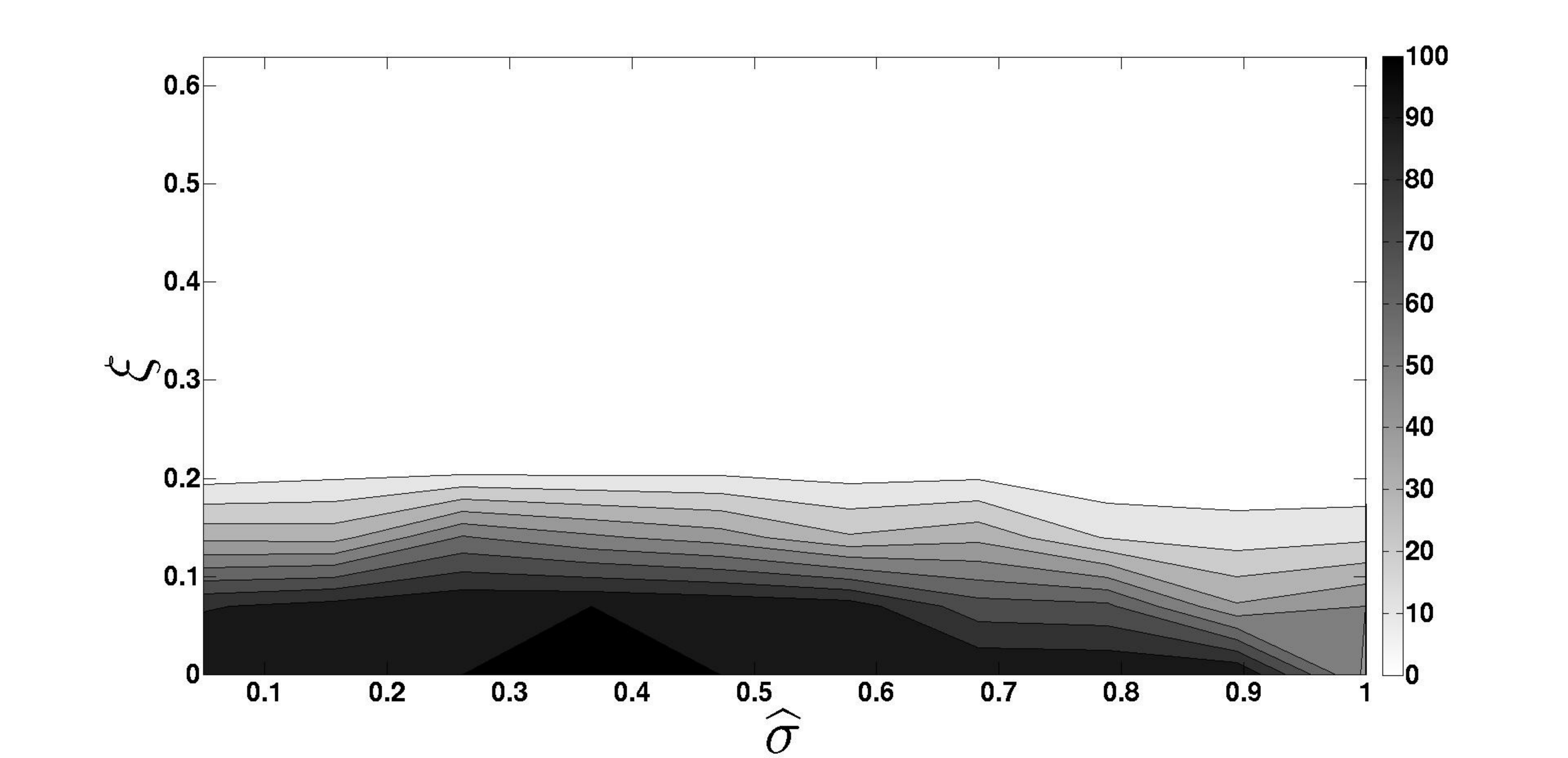}}
}
\subfigure[True parameter values: $\sigma = .75$, $\alpha = .25$]{
\resizebox{0.48\columnwidth}{!}{\includegraphics{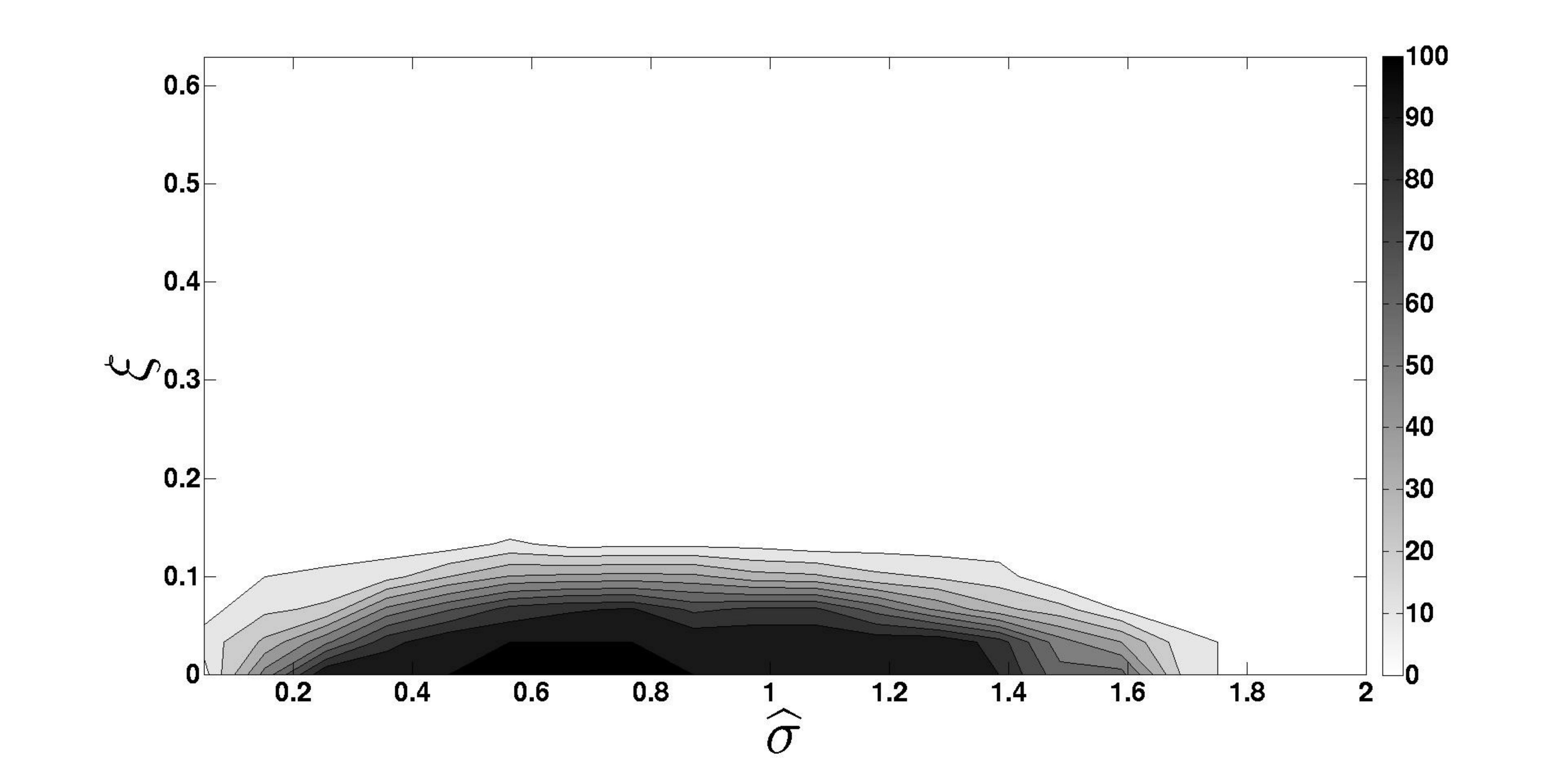}}
}
\caption{Recovery Probabilities when $\alpha = .25$ for two true $\sigma$ settings.  The shade in each phase diagram corresponds to the probability that the greedy algorithm will correctly recover a randomly selected bar code, as a function of the additive noise standard deviation, $\xi$, and the $\sigma$ estimate, $\widehat{\sigma}$.  Black represents correct bar code recovery with probability 1, while pure white represents recovery with probability 0.  Each data point's shade (i.e., probability estimate) is based on 100 random trials.} 
\label{fig:phase2}
\end{figure}

The phase diagrams in Figures~\ref{fig:phase1} and~\ref{fig:phase2}
more clearly illustrate how the reconstruction capabilities of the
greedy algorithm depend on $\sigma$, $\alpha$, the estimate of
$\sigma$, and on the noise level.  We again consider Gaussian additive
noise on the signal, i.e. we consider the inverse problem $d = \alpha
{\cal G}(\sigma){\cal D}x + h$, with independent and identically
distributed $h_j \sim {\cal N}(0,\xi^2)$, for several noise standard
deviation levels $\xi \in [0, .63]$.  Note that
$\mathbb{E}\big(\left\| h |_{I_j} \right\|_1 \big) = 7 r \xi \sqrt{2 /
  \pi}$.\footnote{This follows from the fact that the first raw
  absolute moment of each $h_j$, $\mathbb{E}( | h_j |)$, is $\xi
  \sqrt{2 / \pi}$.}  Thus, the numerical results are consistent with
the bounds in Remark \ref{rek:empirical}.  Each phase diagram
corresponds to different underlying parameter values $(\sigma,
\alpha)$, but in all diagrams we fix the oversampling ratio at $r =
m/n = 10$.  As before, the black regions in the phase diagrams
indicate parameter values $(\widehat{\sigma}, \xi)$ for which $100$
out of $100$ randomly generated bar codes were reconstructed, and
white regions indicate parameter values for which $0$ out of $100$
randomly generated bar codes were reconstructed.

\begin{figure}[]
\subfigure[True parameter values: $\sigma = .45$, $\alpha = 1$.
  Estimated $\widehat{\sigma} = .3$ and Noise Standard Deviation $\xi
  = .3$.  Solving the least-squares problem \eqref{alphals} yields an
  $\alpha$ estimate of $\widehat{\alpha} = .9445$ from
  $\widehat{\sigma}$.  The relative noise-to-signal level, $\nu =
  \frac{\| h \|_2}{\| \alpha {\cal G}(\sigma){\cal D} x \|_2}$, is
  0.4817.]{
  \resizebox{0.98\columnwidth}{!}{
     \includegraphics{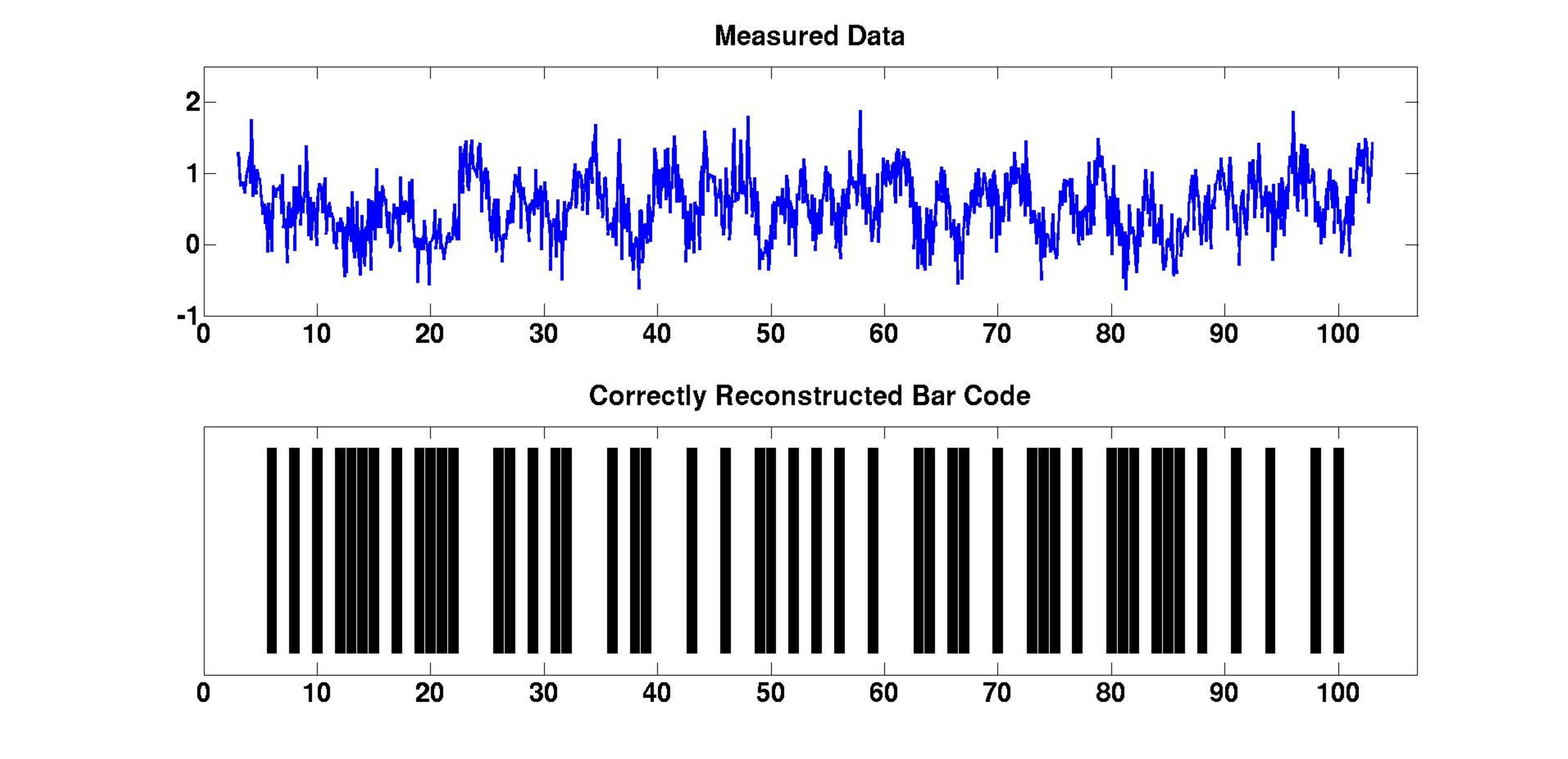}}}

\subfigure[True parameter values: $\sigma = .75$, $\alpha = 1$.
  Estimated $\widehat{\sigma} = 1$ and Noise Standard Deviation $\xi =
  .2$.  Solving the least-squares problem \eqref{alphals} yields an
  $\alpha$ estimate of $\widehat{\alpha} = 1.1409$ from
  $\widehat{\sigma}$.  The relative noise-to-signal level, $\nu =
  \frac{\| h \|_2}{\| \alpha {\cal G}(\sigma){\cal D} x \|_2}$, is
  $0.3362$.]{
  \resizebox{0.98\columnwidth}{!}{
     \includegraphics{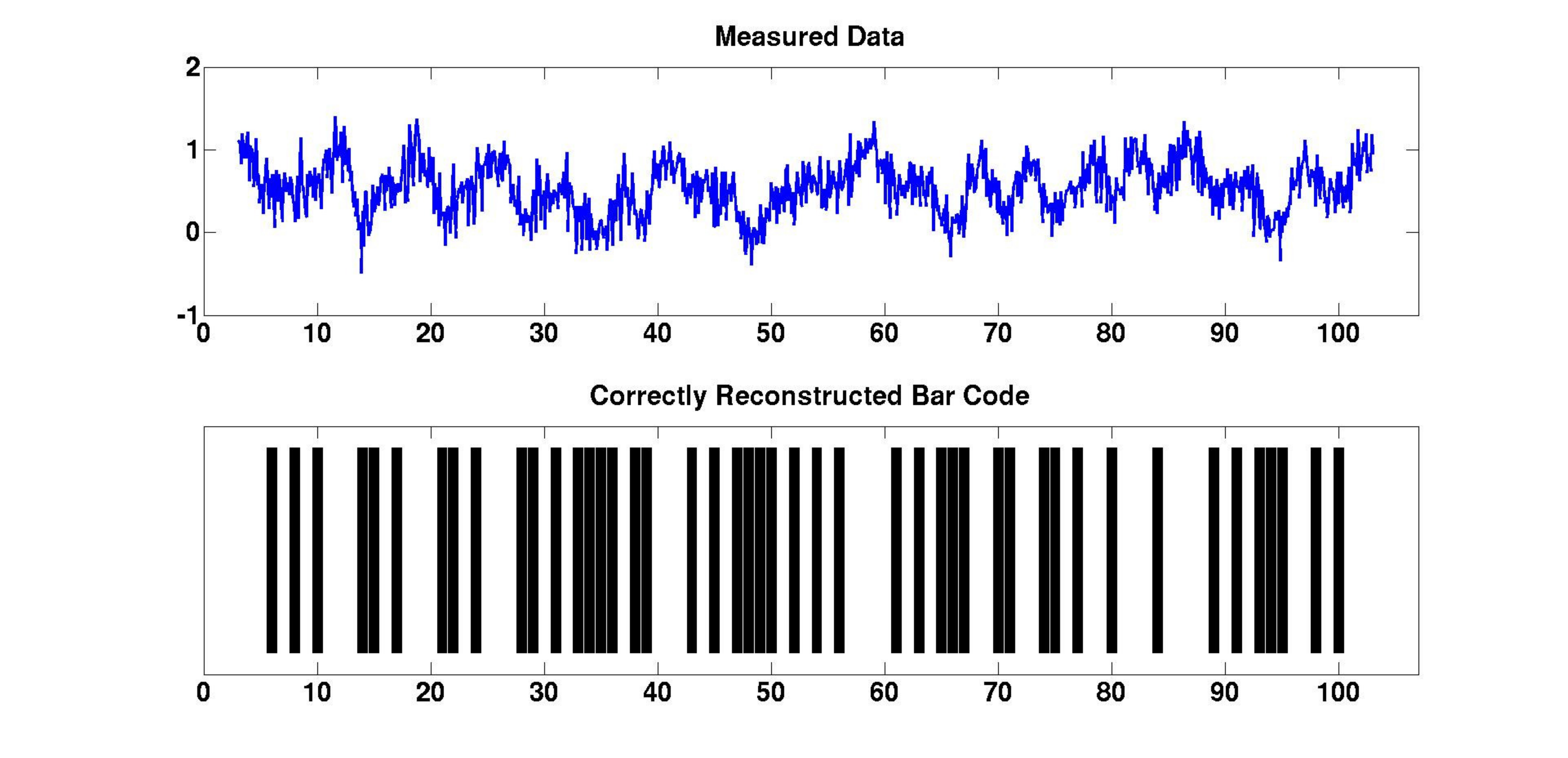}}}
\caption{Two example recovery problems corresponding to dark regions
  in each of the phase diagrams of Figure~\ref{fig:phase1}.  These
  recovery problems are examples of problems with $\alpha = 1$ for
  which the greedy algorithm correctly decodes a randomly selected UPC
  bar code approximately $80\%$ of the time.}
\label{fig:ReconExamples_a1}
\end{figure}

\begin{figure}[]
\subfigure[True parameter values: $\sigma = .45$, $\alpha = .25$.
  Estimated $\widehat{\sigma} = .5$ and Noise Standard Deviation $\xi
  = .1$.  Solving the least-squares problem \eqref{alphals} yields an
  $\alpha$ estimate of $\widehat{\alpha} = 0.2050$ from
  $\widehat{\sigma}$.  The relative noise-to-signal level, $\nu =
  \frac{\| h \|_2}{\| \alpha {\cal G}(\sigma){\cal D} x \|_2}$, is
  0.7001.]{
  \resizebox{0.98\columnwidth}{!}{
     \includegraphics{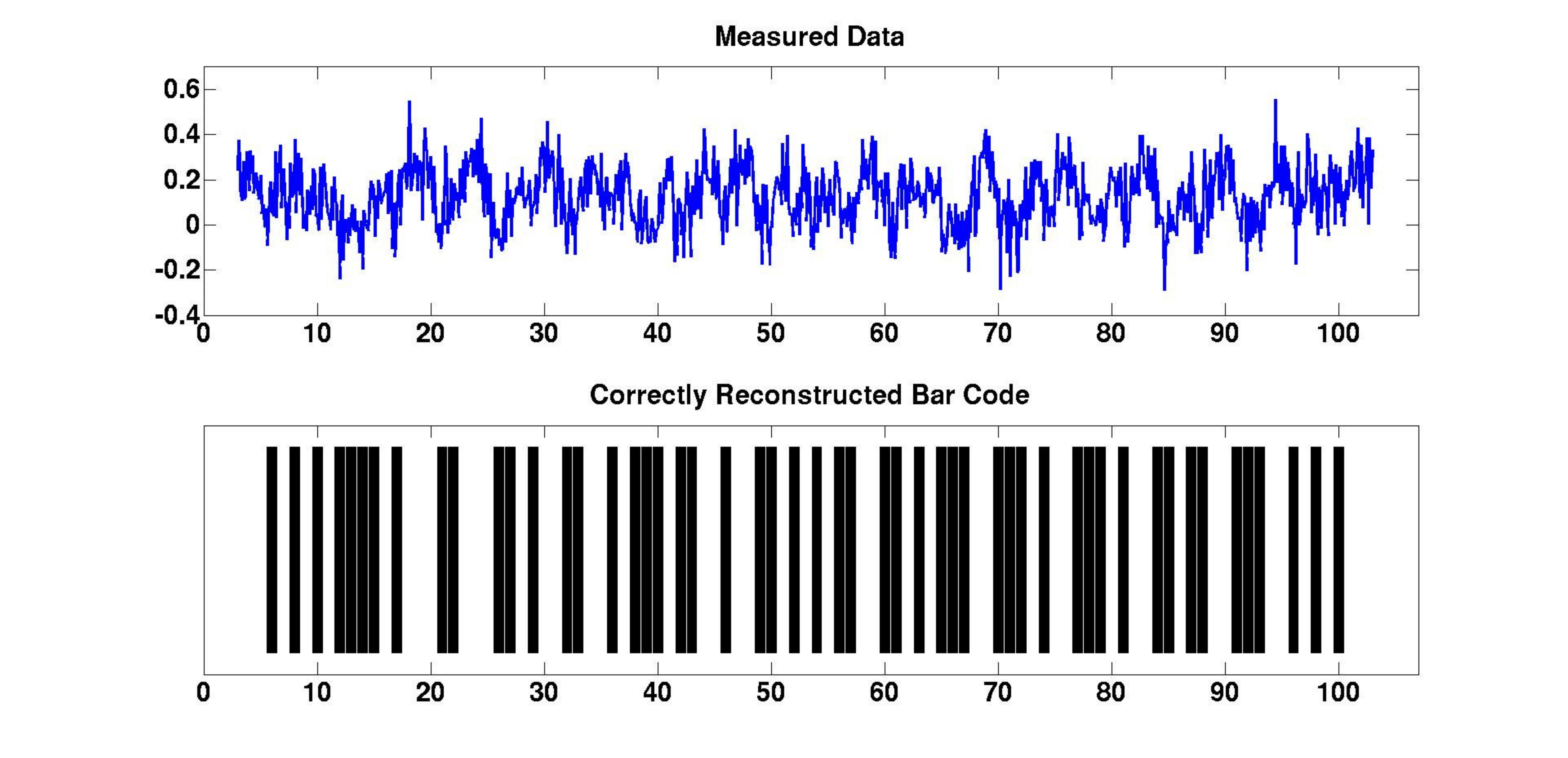}}}

\subfigure[True parameter values: $\sigma = .75$, $\alpha = .25$.
  Estimated $\widehat{\sigma} = .8$ and Noise Standard Deviation $\xi
  = .06$.  Solving the least-squares problem \eqref{alphals} yields an
  $\alpha$ estimate of $\widehat{\alpha} = 0.3057$ from
  $\widehat{\sigma}$.  The relative noise-to-signal level, $\nu =
  \frac{\| h \|_2}{\| \alpha {\cal G}(\sigma){\cal D} x \|_2}$, is
  0.4316.]{
  \resizebox{0.98\columnwidth}{!}{
     \includegraphics{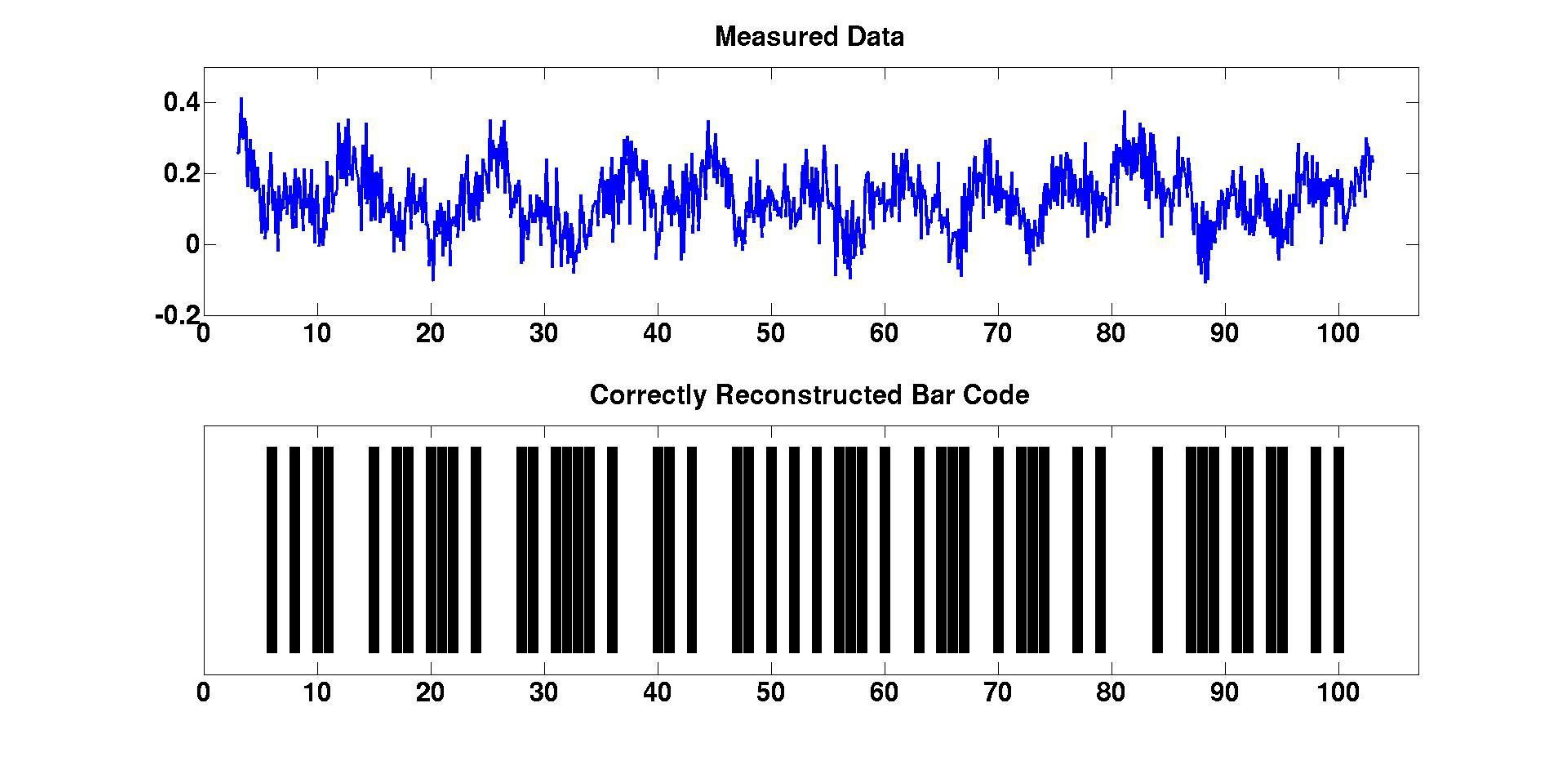}}}

\caption{Two example recovery problems corresponding to dark regions
  in each of the phase diagrams of Figure~\ref{fig:phase2}.  These
  recovery problems are examples of problems with $\alpha = .25$ for
  which the greedy algorithm correctly decodes a randomly-selected UPC
  bar code approximately $60 \%$ of the time.}
\label{fig:ReconExamples_a025}
\end{figure}

Comparing Figures~\ref{fig:phase1}(a) and~\ref{fig:phase2}(a) with
Figures~\ref{fig:phase1}(b) and~\ref{fig:phase2}(b), respectively, we
can see that the greedy algorithm's performance appears to degrade
with increasing $\sigma$.  Note that this is consistent with our
analysis of the algorithm in Section~\ref{sec:AlgAnal}.  Increasing
$\sigma$ makes the forward map $P = \alpha {\cal G}(\sigma){\cal D}$
less block diagonal, thereby increasing the effective value of
$\varepsilon$ in conditions~(\ref{energyconc1})
and~(\ref{energyconc2}).  Hence, condition~(\ref{VictoryCond2}) will be
less likely satisfied as $\sigma$ increases.

Comparing Figures~\ref{fig:phase1} and~\ref{fig:phase2} reveals the
effect of $\alpha$ on the likelihood that the greedy algorithm
correctly decodes a bar code.  As $\alpha$ decreases from $1$ to $.25$
we see a corresponding deterioration of the greedy algorithm's ability
to handle additive noise of a given fixed standard deviation.  This is
entirely expected since $\alpha$ controls the magnitude of the blurred signal $\alpha {\cal G}(\sigma){\cal D} x$.  Hence,
decreasing $\alpha$ effectively decreases the signal-to-noise ratio of
the measured input data $d$.

Finally, all four of the phase diagrams in Figures~\ref{fig:phase1} and~\ref{fig:phase2} demonstrate how the greedy algorithm's probability of successfully recovering a randomly selected bar code varies as a function of the noise standard deviation, $\xi$, and $\sigma$ estimation error, $|\widehat{\sigma}-\sigma|$.  As both the noise level and $\sigma$ estimation error increase, the performance of the greedy algorithm smoothly degrades.  Most importantly, we can see that the greedy algorithm is relatively robust to inaccurate $\sigma$ estimates at low noise levels.  When $\xi \approx 0$ the greedy algorithm appears to suffer only a moderate decline in reconstruction rate even when $|\widehat{\sigma}-\sigma| \approx \sigma$.

Figure~\ref{fig:ReconExamples_a1} gives examples of two bar codes which the greedy algorithm correctly recovers when $\alpha = 1$, one for each value of $\sigma$ presented in Figure~\ref{fig:phase1}.  In each of these examples the noise standard deviation, $\xi$, and estimated $\sigma$ value, $\widehat{\sigma}$, were chosen so that they correspond to dark regions of the example's associated phase diagram in Figure~\ref{fig:phase1}.  Hence, these two examples represent noisy recovery problems for which the greedy algorithm correctly decodes the underlying UPC bar code with relatively high probability.\footnote{The $\xi$ and $\widehat{\sigma}$ values were chosen to correspond to dark regions in a Figure~\ref{fig:phase1} phase diagram, not necessarily to purely black regions.}  Similarly, Figure~\ref{fig:ReconExamples_a025} gives two examples of two bar codes which the greedy algorithm correctly recovered when $\alpha = 0.25$.  Each of these examples has parameters that correspond to a dark region in one of the Figure~\ref{fig:phase2} phase diagrams.

\section{Discussion}

In this work, we present a greedy algorithm for the recovery of bar
codes from signals measured with a laser-based scanner.  So far we
have shown that the method is robust to both additive Gaussian noise
and parameter estimation errors.  There are several issues that we
have not addressed that deserve further investigation.

First, we assumed that the start of the signal is well determined.  By
the start of the signal, we mean the time on the recorded signal that
corresponds to when the laser first strikes a black bar.  This
assumption may be overly optimistic if there is a lot of noise in the
signal.  Preliminary numerical experiments suggest that the algorithm
is not overly sensitive to uncertainties in the start time, and we are
currently working on the development of a fast preprocessing algorithm
for locating the start position from the samples.

Second, while our investigation shows that the algorithm is not
sensitive to the parameter $\sigma$ in the model, we did not
address the best means
for obtaining reasonable approximations to $\sigma$.
In applications where the scanner distance from the bar code may vary
(e.g., with handheld scanners) other techniques for determining
$\widehat{\sigma}$ will be required.  Given the robustness of the
algorithm to parameter estimation errors it may be sufficient to
simply fix $\widehat{\sigma}$ to be the expected optimal $\sigma$
parameter value in such situations.  In situations where more accuracy
is required, the hardware might be called on to provide an estimate of
the scanner distance from the bar code it is scanning, which could
then be used to help produce a reasonable $\widehat{\sigma}$ value.
In any case, we leave more careful consideration of methods for
estimating $\sigma$ to future work.

The final assumption we made was that the intensity distribution
is well modeled by a Gaussian.  This may not be sufficiently accurate
for some distances between the scanner and the bar code.  Since
intensity profile as a function of distance can be measured, one
can conceivably refine the Gaussian model to capture the true
behavior of the intensities.

\end{document}